%% file: BB13.tex
\documentclass[12pt]{article}

\usepackage{mathptmx,amsmath,amssymb,amsfonts,amsthm,a4wide}
\usepackage{graphicx,color,xcolor,pictex,epstopdf,url,algorithm,algorithmic,caption, comment}
\newtheorem{theorem}{Theorem}[section]
\newtheorem{lemma}{Lemma}[section]
\newtheorem{corollary}{Corollary}[section]

\newtheorem{remark}{Remark}[section]
\definecolor{brdarkblue}{rgb}{0.0, 0.0, 0.7}
\definecolor{brdarkgreen}{rgb}{0.0, 0.7, 0.0}
\definecolor{bkdarkcyan}{rgb}{0.0, 0.5, 0.5}
\definecolor{bkcyan}{rgb}{0.0, 0.6, 0.7}

\newcommand{\Margin}[1]{}
\newcommand{\revision}[1]{#1}
\newcommand{\rerevision}[1]{{#1}}

\begin{document}
\title{Determining damping terms in fractional wave equations}
\author{Barbara Kaltenbacher\footnote{
Department of Mathematics,
Alpen-Adria-Universit\"at Klagenfurt.
barbara.kaltenbacher@aau.at.}
\and
William Rundell\footnote{
Department of Mathematics,
Texas A\&M University,
Texas 77843. 
rundell@math.tamu.edu}
}
\date{\vskip-3ex}
\maketitle
\begin{abstract}
This paper deals with the inverse problem of recovering an arbitrary number of
fractional damping terms in a wave equation.
We develop several approaches on uniqueness and reconstruction, some of them
relying on Tauberian theorems that provide relations between the asymptotic
behaviour of solutions in time and Laplace domains.
The possibility of additionally recovering space-dependent coefficients
or initial data is discussed.
The resulting methods for reconstructing coefficients and fractional orders in
these terms are tested numerically. 
In addition, we provide an analysis of the forward problem consisting of
a multiterm fractional wave equation.
\end{abstract}

\leftline{\small \qquad\qquad
{\bf Keywords:} 
Inverse problems, multiterm fractional PDE, fractionally damped wave equation}
\smallskip
\leftline{\small \qquad\qquad 
\textbf{{\textsc msc} 2020} {\bf classification:} 35R30, 35R11, 35L10, 78A46}

\input introduction

\input forward

\input uniqueness

\input recon_meth
\input reconstructions

\subsection*{Acknowledgment}
\vskip-4pt
\noindent
The work of the first author was supported by the Austrian Science Fund {\sc fwf}
under the grants P30054 and {\sc doc}78.

\noindent
The work of the second author was supported
in part by the
National Science Foundation through award {\sc dms}-2111020.


\end{document}

%% file: introduction.tex
\section{Introduction}

Wave phenomena like seismic waves or ultrasound propagation are known to exhibit power law frequency-dependent damping,
which can be modelled by time fractional derivative terms in the corresponding {\sc pde}s.
As a matter of fact, the adequate representation of attenuation with power law 
frequency dependence was perhaps the first application of fractional time 
derivative concepts and among the driving forces for their development, 
see Caputo \cite{Caputo:1967} and Caputo and Mainardi~\cite{CaputoMainardi:1971a}.
 
Fractionally damped wave equations arise by combining classical balance laws
with fractional constitutive relations.
More precisely, consider the equation of motion (resulting from a balance of forces) and the representation of the mechanical strain as the symmetric gradient of the mechanical displacement 
\begin{equation}\label{eqn:motion-strain}
\rho\vec{u}_{tt}=\mbox{div}\sigma+\vec{f}\,, \quad \epsilon=\frac12(\nabla\vec{u}+(\nabla\vec{u})^T)
\end{equation}
where $\rho$ is the mass density, $\vec{u}(x,t)$ the vector of displacements and $\sigma(x,t)$ the stress tensor. The crucial constitutive equation is now the material law that relates stress $\sigma$ and strain $\epsilon$.
Some often encountered instances are 
the fractional Newton (also called Scott-Blair model), Voigt, Maxwell, and Zener models
\begin{equation*}
\sigma = b_1 \partial_t^{\alpha} \epsilon\,, \quad
\sigma = b_0 \epsilon + b_1 \partial_t^{\alpha}\epsilon\,, \quad 
\sigma + d_1\partial_t^{\gamma}\sigma = b_0 \epsilon\,, \quad
\sigma + d_1\partial_t^{\gamma}\sigma = b_0 \epsilon + b_1 \partial_t^{\alpha}\epsilon\,.
\end{equation*}
For a survey on fractionally damped acoustic wave equations we refer to \cite{CaiChenFangHolm_survey2018}.
Driven by applications that are not covered by these models, a general model class has been defined by 
\begin{equation}\label{eqn:sumfrac}
\sum_{j=0}^J d_j \partial_t^{\gamma_j} \sigma = \sum_{k=0}^K b_k \partial_t^{\alpha_k} \epsilon
\end{equation}
with the normalisation $d_0=1$, where $0\leq\gamma_0<\gamma_1\cdots\leq \gamma_J$ and $0\leq\alpha_0<\alpha_1\cdots\leq \alpha_K$
cf., e.g., \cite[section 3.1.1]{Atanackovic_etal:2014}, \cite{bagley1986fractional,schmidt2002finite}.
Combination of \eqref{eqn:sumfrac} with \eqref{eqn:motion-strain} leads to a general wave equation, which  
for simplicity of exposition in the scalar case of one-dimensional mechanics or in acoustics reads as 
\[
\sum_{j=0}^J d_j \partial_t^{2+\gamma_j} u - \sum_{k=0}^K b_k \partial_t^{\alpha_k} \triangle u = r\,.
\]
In this more complicated situation, direct reconstruction of the damping terms from experimental data is not possible any more. It is the aim of this paper to investigate reconstruction of the orders $\alpha_k$, $\gamma_j$ and the coefficients $b_k$, $d_j$ from time trace measurements of $u$, which plays the role of a displacement in mechanics or of a pressure in acoustics.
In the latter application case, on which we mainly focus here, certain imaging modalities additionally require  reconstruction of an unknown sound speed $c(x)$ in $b_0=c^2$ (ultrasound tomography), and/or an initial condition $u_0(x)$ (equivalently, of the spatially varying part $f$ of the source $r(x,t)=\sigma(t)f(x)$; photoacoustic tomography).
Besides investigating these tasks, inspired by recent work by Jin and Kian \cite{JinKian:2021} in the context of diffusion models, we will also study the question of identifiability of the damping terms in an unknown medium, that is, e.g., with unknown $c(x)$, without aiming at its reconstruction. 
While to the best of our knowledge, this is the first work on recovering
multiple fractional damping in wave type equations, much more literature is
available for anomalous diffusion {\sc pde}s.
Besides the reference \cite{JinKian:2021},
we also refer to, e.g.,
\cite{HatanoNakagawaWangYamamoto:2013,LiYamamoto:2015,LiZhangJiaYamamoto:2013,Yamamoto:2021,SunLiZhang:2021}.

We will present results on uniqueness and on reconstruction algorithms in
each of three paradigms depending on the availability of measured data.
These are:  a full time trace $u(x_0,t)$, $t >0$ for some fixed point $x_0$;
extreme large time measurements; and extremely small time measurements.
The main tool used in the latter two cases is a Tauberian theorem that
relates the asymptotic behaviour of the solution in time to the structure
of the representation of the solution in the Laplace transformed formulation
that contains specific algebraic information on the fractional operator.

The remainder of this paper is organised as follows.
In subsections~\ref{sec:models}, \ref{sec:inverse} we introduce the models under consideration and state the inverse problem.
Section~\ref{sec:forward} is devoted to the forward problem: We provide an analysis of the general midterm fractional wave equation and describe the numerical forward solver used in our reconstructions.
In section~\ref{sec:uniqueness} we prove uniqueness results based on different paradigms: Using the Weyl estimate on decay of eigenvalues together with smoothness of the excitation, applying single mode excitations and exploiting information on poles and residues in the Laplace domain.
Section~\ref{sec:recon_meth} develops reconstruction methods in the three scenarios indicated above: full time (that is, Laplace domain), large time and small time observations. 
Correspondingly, in section~\ref{sec:reconstructions} we provide numerical reconstructions along with a discussion of the approaches.

\input models
\input inverse

%% file: models.tex
\subsection{Models} \label{sec:models}
We consider the second order in time damped wave equation
\begin{equation}\label{eqn:general_2ndorder}
u_{tt}+c^2\mathcal{A}u+\sum_{k=1}^N b_k \partial_t^{\alpha_k} \mathcal{A}^{\beta_k}u=r
\end{equation}
where $\mathcal{A}=-\mathbb{L}$ on $\Omega\subseteq\mathbb{R}^d$, equipped with homogeneous Dirichlet/ Neumann/ impedance boundary conditions, for an elliptic differential operator $\mathbb{L}$ and 
$\alpha_k\in(0,1]$, $\beta_{k}\in(\frac12,1]$.

Throughout this paper, we denote by $\partial_t^\alpha$ the (partial) Caputo-Djrbashian fractional time derivative of order $\alpha\in(n-1,n)$ with $n\in\mathbb{N}$ by
\[
\partial_t^\alpha u  = \frac{1}{\Gamma(n-\alpha)}\int_{0}^{t}\frac{\partial^n_t u(s)}{(t-s)^{\alpha+1-n}}  \,ds\;,
\]
where $\partial^n_t$ denotes the $n$-th integer order partial time derivative and for
$\gamma\in(0,1)$, and $I_t^\gamma$ is the Abel integral operator defined by 
$
I_t^\gamma[v](t) = \frac{1}{\Gamma(\gamma)}\int_{0}^{t}\frac{v(s)}{(t-s)^{1-\gamma}}ds\,.
$
For details on fractional differentiation and subdiffusion equations,
we refer to, e.g.,
\cite{Dzjbashian:1966,Djrbashian:1993,
MainardiGorenflo:2000,SakamotoYamamoto:2011a,SamkoKilbasMarichev:1993}.
See also the tutorial paper on inverse problems for anomalous diffusion processes
\cite{JinRundell:2015}.

We assume that the differentiation orders with respect to time are distinct, that is
\begin{equation}\label{eqn:alphas_distinct}
0<\alpha_1<\alpha_2<\cdots<\alpha_N\leq1,\quad
0<\gamma_1<\gamma_2<\cdots<\gamma_J\leq1
\end{equation}
a property that is crucial for distinguishing the different asymptotic terms in the solution $u$ and its Laplace transform.

More generally than \eqref{eqn:general_2ndorder},
in most of this paper we also take higher than second order time derivatives into account
\begin{equation}\label{eqn:general}
u_{tt}+c^2\mathcal{A}u+\sum_{j=1}^J d_j\partial_t^{2+\gamma_j}u+\sum_{k=1}^N b_k \partial_t^{\alpha_k} \mathcal{A}^{\beta_k}u=r
\end{equation}

Note that the usual models from \cite{KaltenbacherRundell:2021b} are contained in \eqref{eqn:general}, which is actually an extension of the general model from, e.g., the book \cite{Atanackovic_etal:2014}
in the sense that it also allows for fractional powers of the operator $\mathcal{A}$.

Typically, in acoustics we just have $\mathbb{L}=\triangle$, and then the operator $\mathcal{A}$ is known. To take into account a (possibly unknown) spatially varying speed of sound $c(x)$, one can instead consider 
\begin{equation}\label{eqn:general_c}
u_{tt}+\mathcal{A}_c u+\sum_{j=1}^J d_j\partial_t^{2+\gamma_j}u+\sum_{k=1}^N b_k \partial_t^{\alpha_k} \mathcal{A}_c^{\beta_k}u=r
\end{equation}
with $\mathcal{A}_c=-c(x)^2\triangle$; in order to get a selfadjoint operator $\mathcal{A}_c$ we then use the weighted $L^2$ inner product with weight function $\frac{1}{c^2}$.

\medskip

Note that in order to assume \eqref{eqn:alphas_distinct} without loss of generality, so allowing for all possible combinations of $\alpha_k$, $\beta_k$, we would have to consider
\begin{equation}\label{eqn:mostgeneral_2ndorder}
u_{tt}+c^2\mathcal{A}u+\sum_{k=1}^N\partial_t^{\alpha_k}\sum_{\ell=1}^{M_k} c_{k\ell}\mathcal{A}^{\beta_{k\ell}} u=r
\end{equation}
or, including higher than second time derivatives such as in \eqref{eqn:general}
\begin{equation}\label{eqn:mostgeneral}
u_{tt}+c^2\mathcal{A}u+\sum_{j=1}^J d_j\partial_t^{2+\gamma_j}u+\sum_{k=1}^N\partial_t^{\alpha_k}\sum_{\ell=1}^{M_k} c_{k\ell}\mathcal{A}^{\beta_{k\ell}}u=r.
\end{equation}
However, \eqref{eqn:mostgeneral_2ndorder}, \eqref{eqn:mostgeneral} might be
viewed as over-parameterised models and they actually lead to difficulties
in proving uniqueness of solutions.


\medskip

Important special cases of \eqref{eqn:general_c} are the Caputo-Wismer-Kelvin-Chen-Holm 
\begin{equation}\label{eqn:CH_c}
u_{tt}+\mathcal{A}_c u+ b \partial_t^{\alpha} \mathcal{A}_c^\beta u=r
\end{equation}
and the fractional Zener {\sc fz} model
\begin{equation}\label{eqn:FZ_c}
u_{tt}+\mathcal{A}_c u+ d \partial_t^{2+\gamma}  u + b \partial_t^{\alpha} \mathcal{A}_c u=r .
\end{equation}

These {\sc pde} models will be considered on a time interval $t\in(0,T)$ and driven by initial conditions and/or a separable source term
\begin{equation}\label{eqn:init_source}
\begin{aligned}
&r(x,t)=\sigma(t)f(x)\,, \quad x\in\Omega, \quad t\in(0,T)\\
&u(x,0)=u_0(x)\,, \quad u_t(x,0)=u_1(x)\,, \quad (u_{tt}(x,0)=u_2(x)\mbox{ if }\gamma_J>0)\quad x\in\Omega. \end{aligned}
\end{equation}
Note that boundary conditions are already incorporated into the operator $\mathcal{A}$ or $\mathcal{A}_c$, respectively. 

%% file: inverse.tex
\subsection{Inverse problems} \label{sec:inverse}
\rerevision{
As in \cite{JinKian:2021}, we assume that data has been obtained from a quite general experimental setup, in which the {\sc pde} \eqref{eqn:general} with its unknown coefficients remains the same, while both excitation \eqref{eqn:init_source} and observation location may vary between the experiments. We number the experiments with an index $i$ to denote the corresponding excitation by $(u_{0,i},u_{1,i},f_i)$, and the corresponding observation operator by $B_i$. Given $I$ experiments, we thus 
assume to have the following $I$ time trace observations for different driving by initial and/or source data:
\begin{equation}\label{eqn:obs_titr}
\begin{aligned}
&h_i(t)=(B_i u_i)(t), \ t\in(0,T) \quad\mbox{ where $u_i$ solves \eqref{eqn:general}, \eqref{eqn:init_source} with } (u_0,u_1,f)=(u_{0,i},u_{1,i},f_i)\,,\\
&i=1,\ldots I.
\end{aligned}
\end{equation}
Examples of observation operators are 
\[
(a) \ \; B_i v= v(x_i)\quad \mbox{ or }\quad 
(b) \ \; B_i v= \int_{\Sigma_i}\eta_i(x) v(x)\, dx
\]
for some points $x_i\in\overline{\Omega}$ and some weight functions
$\eta_i\in L^\infty(\omega)$, $\Sigma_i\subseteq\partial\Omega$, provided that these evaluations are well-defined, that is, $u(t)\in C(\overline{\Omega})$ in case (a), or $u(t)\in L^1(\omega)$ in case (b), which according to Theorem~\ref{thm:sumfrac} is the case if $\Omega\subseteq\mathbb{R}^d$ with $d=1$ for case (a) or $d\in\{1,2,3\}$ for case (b). 
}
\Margin{Ref 2 (i)}

The coefficients in \eqref{eqn:general} that we aim to in recover are 
\[
N, \ J, \ c, \ b_k, \ d_j, \ \alpha_k, \ \beta_k, \ \gamma_j, \quad
j=1,\ldots,J, \ k=1,\ldots,N\;.
\]

Moreover, we consider the problem of identifying, besides the differentiation orders, some space dependent quantities. 
In practice the most relevant ones are either the speed of sound $c=c(x)$ in ultrasound tonography, or the initial data $u_0=u_0(x)$ while $u_1=0$ (equivalently the source term $f=f(x)$, while $u_0=0$, $u_1=0$) in photoacoustic tomography {\sc pat}.
For the latter purpose, we will consider observations not only at single points (a) or patches (b), but over a surface $\Sigma$.

%% file: forward.tex
\section{The forward problem}\label{sec:forward}

\input analysis_forward

\input numerics_forward

\input resolvent

%% file: analysis_forward.tex
\def\calAtil{\tilde{\mathcal{A}}}

\subsection{Analysis of the forward problem}

In this section we prove well-posedness of the forward problem \eqref{eqn:general_c}.
Just as in the inverse problem which is the main topic of this paper,
we focus on the case of constant coefficients $b_k$, $d_j$, while the sound speed $c=c(x)$ contained in the operator $\calAtil$ may vary in space. We refer to Chapter 7 of the upcoming book \cite{BBB} for an analysis in case all coefficients are space dependent.

Our analysis will be based on energy estimates obtained by multiplying the {\sc pde} (or actually its Galerkin semidiscretization with respect to space) with appropriate multipliers. For this purpose, we will make use of certain mapping properties of fractional derivative operators that can be quantified via upper and lower bounds.

We will use the following results that hold for $\theta\in[0,1)$.
\begin{itemize}	
	\item
	From \cite[Lemma 2.3]{Eggermont:1987}; see also \cite[Theorem 1]{VoegeliNedaiaslSauter:2016}: For any $w\in H^{-\theta/2}(0,t)$,
	\begin{equation}\label{eqn:coercivityI}
	\int_0^t \langle {_0I_t}^{\theta} w(s),  w(s) \rangle ds \geq \cos ( \pi\theta/2 ) \| w \|_{H^{-\theta/2}(0,t)}^2, 
	\end{equation}	 
	\item 
	From \cite[Theorems 2.1, 2.2, 2.4 and Proposition 2.1]{KubicaRyszewskaYamamoto:2020}:
	There exist constants $0<\underline{C}(\theta)<\overline{C}(\theta)$ such that for any $w\in H^\theta(0,t)$ (with $w(0)=0$ if $\theta\geq\frac12$), the equivalence estimates 
	\begin{equation}\label{eqn:equivalence_pos}
\underline{C}(\theta) \|w\|_{H^\theta(0,t)} \leq \|\partial_t^\theta w\|_{L^2(0,t)}\leq
\overline{C}(\theta) \|w\|_{H^\theta(0,t)}
	\end{equation}
hold.
Likewise, for any $w\in L^2(0,t)$ we have $\partial_t^\theta w\in H^{-\theta}(0,t)$ and 
	\begin{equation}\label{eqn:equivalence_neg}
\underline{C}(\theta) \|w\|_{L^2(0,T)} \leq \|\partial_t^\theta w\|_{H^{-\theta}(0,t)}\leq
\overline{C}(\theta) \|w\|_{L^2(0,t)}.
	\end{equation}
Moreover, (by a duality argument) there exists a constant $C(2\theta)>0$ such that for any $w\in L^2(0,t)$
	\begin{equation}\label{eqn:boundednessI}
	\int_0^t ({_0I_t}^{2\theta} w(s))  w(s) ds \leq C(2\theta)\| w \|_{H^{-\theta}(0,t)}^2, 
	\end{equation}	 
\end{itemize}
\rerevision{Note that the latter two facts \eqref{eqn:equivalence_neg}, \eqref{eqn:boundednessI} hold true for $\theta=1$ as well.}

We are now prepared to state and prove our main theorem on well-posedness of an initial value problem 
\begin{equation}\label{eqn:sumfrac_wave_acou_ibvp}
  \left\{\begin{aligned}
\sum_{j=0}^J d_j &\partial_t^{2+\gamma_j} u + \sum_{k=0}^K b_k \partial_t^{\alpha_k} 
\calAtil u = r\quad \mbox{ in }\Omega\times(0,T),\\
    &u(x,0)  = u_0(x)\,, \ u_t(x,0) = u_1(x) \quad x\in\Omega,\\
&(\mbox{if }\gamma_J>0:\ u_{tt}(x,0) = u_2(x)\quad x\in\Omega,)
  \end{aligned}\right.
\end{equation}
which clearly covers the models from Section~\ref{sec:models} in case of $\beta_1=\cdots=\beta_K$.
Here $\calAtil$ is a selfadjoint (with respect to some possibly weighted space $L^2_w(\Omega)$ with inner product $\langle\cdot,\cdot\rangle$) operator with
eigenvalues and eigenfunctions $(\lambda_i,\varphi_i)_{i\in\mathbb{N}}$,
and $\dot{H}^1(\Omega)=\{v\in L^2_w(\Omega)\, : \, \|\calAtil^{1/2}v\|_{L^2_w(\Omega)}<\infty\}$.
This includes the cases $\calAtil=-\triangle$, $\calAtil=-c(x)^2\triangle$, $\calAtil=-c(x)^2\nabla\cdot(\tfrac{1}{\rho(x)}\nabla)$ relevant for acoustics, where in the latter two cases we use the inner product $\langle v_1,v_2\rangle=\int_\Omega\frac{1}{c^2(x)}\, v_1(x)v_2(x)\ dx$ in which $\calAtil$ is selfadjoint.

The coefficients and orders will be assumed to satisfy 
\begin{equation}\label{eqn:gammaalpha}
0\leq\gamma_0<\gamma_1<\cdots< \gamma_J
 \leq1
\,, \qquad 0\leq\alpha_0<\alpha_1<\cdots< \alpha_K
 \leq1
\,, 
\qquad \gamma_J\leq\alpha_K\,,
\end{equation}
\Margin{Ref 2 (ii)}
and 
\begin{equation}\label{eqn:dJbK}
d_j\geq0\,, \ j\in\{1,\ldots,J\}\,, \quad b_k\geq 0\,, \ k\in\{\rerevision{k_*},\ldots,K\}\,, \quad d_J>0 \,, \quad b_K>0 \,.
\end{equation} 
with $k_*$ to be specified below.
\rerevision{The condition $\gamma_J\leq\alpha_K$ has been shown to be necessary for thermodynamic consistency, see \cite[Lemma 3.1]{Atanackovic_etal:2014}.\\
Note that the proof of Theorem~\ref{thm:sumfrac} also goes through without the restrictions $\gamma_j \leq1$, $\alpha_k \leq1$, but orders larger than one would require introducing further initial conditions, which would make the exposition less transparent. Moreover, the focus in this paper is actually on wave type equations where this restriction of the orders is physically relevant.  
}

We first of all consider the following time integrated and therefore weaker form with $\tilde{r}(x,t)=\int_0^t r(x,s)\, ds$
\begin{equation}\label{eqn:sumfrac_wave_acou_ibvp_timeint}
  \left\{\begin{aligned}
&\sum_{j=0}^J d_j
\Bigl(\partial_t^{1+\gamma_j} u-\frac{t^{1-\gamma_j}}{\Gamma(2-\gamma_j)}u_2(x)\Bigr)\\
&+ \sum_{k=0}^K b_k \Bigl({_0I_t}^{1-\alpha_k} \calAtil u - \frac{t^{1-\alpha_k}}{\Gamma(2-\alpha_k)}\calAtil u_0\Bigr)
= \tilde{r} \quad \mbox{ in }\Omega\times(0,T),\\
    &\hspace*{1cm}u(x,0)  = u_0(x)\,, \ u_t(x,0) = u_1(x) \quad x\in\Omega
  \end{aligned}\right.
\end{equation}

The results are obtained by testing the {\sc pde} in 
\eqref{eqn:sumfrac_wave_acou_ibvp} 
with 
\rerevision{
$\partial_t^{\tau} u$ for some $\tau>0$
that needs to be well chosen in order to exploit the estimates \eqref{eqn:coercivityI}-\eqref{eqn:boundednessI}.
Generalizing the common choice $\tau=1$ in case of the classical second order wave equation $\gamma_0=\cdots=\gamma_J=\alpha_0=\cdots=\alpha_K=0$ and in view of the stability condition $\gamma_J\leq\alpha_K$ we choose $\tau\in [1+\gamma_J,1+\alpha_K]$, actually, 
$\tau=1+\alpha_{k_*}$ with $k_*$ according to 
\begin{equation}\label{eqn:mast}
k_*\in\{1,\ldots,K\}\mbox{ such that }\alpha_{k_*-1}<\gamma_J\leq\alpha_{k_*},
\end{equation}
if $\alpha_0<\gamma_J$ and $k_*=0$ otherwise.
More precisely, we will assume that 
\begin{equation}\label{eqn:i-or-ii}
(i) \   \alpha_0<\gamma_J, \ k_*\mbox{ as in \eqref{eqn:mast} }, 
\ \mbox{ or } \
(ii)\  
\alpha_0\geq\gamma_J, \ k_*:=0
\end{equation}
holds.
Additionally, we assume
\begin{equation}\label{alphakstar}
\alpha_K-1\leq\alpha_{k_*}\leq 1+\min\{\gamma_0,\alpha_0\}
\ \mbox{ and either } \ \gamma_J<\alpha_K \mbox{ or }\ \gamma_J=\alpha_K=\alpha_0.
\end{equation}
The upper and lower bounds on $\alpha_{k_*}$ in \eqref{alphakstar} are needed to guarantee $\theta\in[0,1)$ in several instances of $\theta$ used in \eqref{eqn:coercivityI}-\eqref{eqn:boundednessI} in the proof below.
\\
Note that we are considering the case $\gamma_J=\alpha_K$ only in the specific setting of a single $\calAtil$ term.
Indeed, when admitting several $\calAtil$ terms, well-posedness depends on specific coefficient constellations. As an example of this, let us mention the fractional Moore-Gibson-Thompson equations arising in acoustics
\begin{equation}\label{fMGT}
u_{tt}+\mathfrak{T}\partial_t^{2+\alpha} u-\mathfrak{c}^2\Delta u-\mathfrak{b}\partial_t^{\alpha} \Delta u = f
\end{equation}
with positive parameters $\mathfrak{T},\mathfrak{b},\mathfrak{c},\alpha$.
It is well-posed only if $\mathfrak{b}\geq\mathfrak{c}^2\mathfrak{T}$, according to \cite[Proposition 7.1]{fracJMGT}.
As the analysis of this example also shows, even lower order time derivatives of $\calAtil u$ (such as the $-\mathfrak{c}^2\Delta u$ term in \eqref{fMGT}) may prohibit global in time boundedness of solutions. 
It will turn out in the proof that we therefore have to impose a bound on the magnitude of the lower order term coefficients in terms of the highest order one
\begin{equation}\label{smallnessbk}
\begin{aligned}
&\sum_{k=0}^{k_*-1} C(1+\alpha_{k_*}-\alpha_k)\overline{C}((1+\alpha_{k_*}-\alpha_k)/2)\, |b_k|\, \|_0I_t^{(\alpha_K-\alpha_k)/2}\|_{L^2(0,T)\to L^2(0,T)}\\
&\leq
\underline{C}((1+\alpha_{k_*}-\alpha_K)/2)\, b_K\, \cos ( \pi(1+\alpha_{k_*}-\alpha_K)/2 ).
\end{aligned}
\end{equation}
that due to $T$ dependence of $\|_0I_t^{(\alpha_K-\alpha_k)/2}\|_{L^2(0,T)\to L^2(0,T)}$ clearly restricts the maximal time $T$ unless all $b_k$ with index lower than $k_*$ vanish, cf. Corollary~\ref{cor:global} below.
}

\rerevision{By applying the test functions as mentioned above}, we will derive estimates for the energy functionals (which are not supposed to represent physical energies but are just used for mathematical purposes)
\begin{equation}\label{eqn:Ea}
\begin{aligned}
\mathcal{E}_\gamma[u](t)=&
\sum_{j=0}^{J-1} \cos ( \pi(1+\gamma_j-\alpha_{k_*})/2 ) \| \sqrt{d_j}\partial_t^{(3+\gamma_j+\alpha_{k_*})/2} u \|_{L^2(0,t;L^2_w(\Omega))}^2 + \underline{\mathcal{E}}_\gamma[u](t)
\end{aligned}
\end{equation}
where 
\begin{equation}\label{eqn:ulEa}
\underline{\mathcal{E}}_\gamma[u](t)=
\begin{cases}\tilde{d} \|\partial_t^{(3+\gamma_J+\alpha_{k_*})/2} u \|_{L^2(0,t;L^2_w(\Omega))}^2
&\mbox{ if }\gamma_J<\alpha_{k_*}\\
\frac{d_J}{4}\|\partial_t^{1+\gamma_J} u \|_{L^\infty(0,t;L^2_w(\Omega))}^2 &\mbox{ if }\gamma_J=\alpha_{k_*}
\end{cases}
\end{equation}
where 
$\tilde{d}=d_J\,\cos ( \pi(1+\gamma_J-\alpha_{k_*})/2 )  \underline{C}((1+\gamma_J-\alpha_{k_*})/2)^2$ with $\underline{C}$ as in \eqref{eqn:equivalence_neg}, 
\rerevision{$\theta=(1+\gamma_J-\alpha_{k_*})/2\in[0,1)$}
and 
\begin{equation}\label{eqn:Eb}
\begin{aligned}
\mathcal{E}_\alpha[u](t)=&
\frac14\|\sqrt{b_{k_*}}\partial_t^{\alpha_{k_*}}\calAtil^{1/2}u\|_{L^\infty(0,t);L^2_w(\Omega)}^2\\
&+\sum_{k=k_*+1}^K
\cos ( \pi(1+\alpha_{k_*}-\alpha_k)/2 ) \| \sqrt{b_k}\partial_t^{(1+\alpha_{k_*}+\alpha_k)/2} \calAtil^{1/2}u \|_{L^2(0,t;L^2_w(\Omega))}^2
\end{aligned}
\end{equation}
where the sum in \eqref{eqn:Eb} is void in case $k_*=K$.

The solution space induced by these energies is $U_\gamma\cap U_\alpha$ where
\begin{equation}\label{eqn:defU}
\begin{aligned}
U_\gamma&=\begin{cases}
H^{(3+\gamma_J+\alpha_{k_*})/2}(0,T;L^2_w(\Omega))&\mbox{ if }\gamma_J<\alpha_{k_*}\\
\{v\in L^2(0,T;L^2_w(\Omega))\, : \, \partial^{1+\gamma_J} v\in L^\infty(0,T;L^2_w(\Omega))&\mbox{ if }\gamma_J=\alpha_{k_*}
\end{cases}
\\
U_\alpha&=\begin{cases}
H^{(1+\alpha_K+\alpha_{k_*})/2}(0,T;\dot{H}(\Omega))&\mbox{ if }K>k_*\\
\{v\in \rerevision{H^{(1-\epsilon)/2+\alpha_K}(0,T;\dot{H}(\Omega))}\, : \, \partial^{\alpha_K} v\in L^\infty(0,T;\dot{H}(\Omega))&\mbox{ if }K=k_*
\end{cases}
\end{aligned}
\end{equation}
\rerevision{for any $\epsilon\in(0,\alpha_K-\max\{\alpha_{K-1},\gamma_J\})$.}

\begin{theorem}\label{thm:sumfrac}
Assume that the coefficients $\gamma_j$, $\alpha_k$, $d_j$, $b_k$ satisfy \eqref{eqn:gammaalpha}, \eqref{eqn:dJbK}, \eqref{eqn:i-or-ii}, \eqref{alphakstar}, 
\rerevision{and that $T$ and/or the coefficients $b_k$ for $k\leq k_*$ are small enough so that \eqref{smallnessbk} holds.}

Then for any $u_0,u_1\in \dot{H}^1(\Omega)$, ($u_2\in L^2_w(\Omega)$ if $\gamma_J>0$), $r\in H^{\alpha_{k_*}-\gamma_J}(0,T;L^2_w(\Omega))$, the time integrated initial boundary value problem
\eqref{eqn:sumfrac_wave_acou_ibvp_timeint} (to be understood in a weak $\dot{H}^{-1}(\Omega)$ sense with respect to space and an $L^2(0,T)$ sense with respect to time) has a unique solution $u\in U_\gamma\cap U_\alpha$ cf. \eqref{eqn:defU}.
This solution satisfies the energy estimate
\begin{equation}\label{eqn:enest_sumfrac}
\begin{aligned}
\mathcal{E}_\gamma[u](t)+\mathcal{E}_\alpha[u](t)
\leq 
C(t) \bigl(\|u_0\|_{\dot{H}^1(\Omega)}^2 + C_0+\|r\|_{H^{\alpha_{k_*}-\gamma_J}(0,t,L^2_w(\Omega))}^2\bigr)\,, \quad t\in(0,T)
\end{aligned}
\end{equation}
for some constant $C(t)$ depending on time, \rerevision{which is bounded as $t\to0$ but} in general grows exponentially as
$t\to\infty$, $\mathcal{E}_\gamma$, $\mathcal{E}_\alpha$ defined as in \eqref{eqn:Ea}, \eqref{eqn:Eb}, and 
\begin{equation}\label{eqn:C0}
C_0 = \begin{cases} 
\|u_1\|_{\dot{H}^1(\Omega)}^2 + \|u_2\|_{L^2_w(\Omega)}^2 &\mbox{ if }\alpha_{k_*}>0\\ 
\|u_1\|_{L^2_w(\Omega)}^2 &\mbox{ if }\alpha_{k_*}=0.
\end{cases} 
\end{equation}
\end{theorem}

\begin{proof}
We exclude the case $\alpha_K=0$ since this via \eqref{eqn:gammaalpha} implies that also $\gamma_J=0$ and we would therefore deal with the conventional second order wave equation whose analysis can be found in textbooks, cf. e.g., \cite{Evans:2010}. Thus for the remainder of this proof we assume $\alpha_K>0$ to hold.

\noindent
{\bf Step 1.} Galerkin discretisation.\\
In order to prove existence and uniqueness of solutions to \eqref{eqn:sumfrac_wave_acou_ibvp},
we apply the usual Faedo-Galerkin approach of discretisation in space with eigenfunctions $\varphi_i$ corresponding to eigenvalues $\lambda_i$ of $\calAtil$, 
$u(x,t)\approx u^L(x,t)=\sum_{i=1}^L u^L_i(t)\varphi_i(x)$ and testing with $\varphi_j$, that is,
\begin{equation}\label{eqn:Galerkin}
\begin{aligned}
&\begin{cases}
\langle \sum_{j=0}^J d_j \partial_t^{2+\gamma_j} u^L(t) + \sum_{k=0}^K b_k \partial_t^{\alpha_k} 
\calAtil u^L(t) - r(t),v\rangle_{L^2_w(\Omega)} = 0 \quad t\in(0,T)\\
\langle u^L(0)-u_0,v\rangle = \langle u^L_t(0)-u_1,v\rangle = 0\,, \quad
(\mbox{if }\gamma_J>0:\ \langle u_{tt} - u_2,v\rangle = 0)
\end{cases}\\
&\mbox{ for all }v\in\mbox{span}(\varphi_1,\ldots, \varphi_L)\,.
\end{aligned}
\end{equation}
To prove existence of a solution $\underline{u}^L\in H^{2+\gamma_J}(0,T;\mathbb{R}^L)$
to \eqref{eqn:Galerkin}, 
we rewrite it as a system of Volterra integral equations for 
\[
\underline{\xi}^L:= \partial_t^{2+\gamma_J}\underline{u}^L = 
\begin{cases}
{_0I_t}^{1-\gamma_J}\underline{u}_{ttt}^L &\mbox{ if }\gamma_J>0\\
\underline{u}_{tt}^L &\mbox{ if }\gamma_J=0\,.
\end{cases}
\]
Unique solvability of the resulting system in $L^2(0,T^L)$ follows from~\cite[Theorem 4.2, p. 241 in \S 9]{GLS90}. Then from	
\[ \left \{
	\begin{aligned}
	\partial_t^{\gamma_J}\underline{u}_{tt}^L =\underline{\xi} \in L^2(0,T^L), \quad
	\underline{u}^L_{tt}(0)=\underline{u}_{2}^L
	\end{aligned} \right.
	\]
in case $\gamma_J>0$ we have a unique $\underline{u}_{tt}^L \in H^{\gamma_J}(0,T^L)$ (cf.~\cite[\S 3.3]{KubicaRyszewskaYamamoto:2020}); the same trivially holds true in case $\gamma_J=0$. 

From the energy estimates below, this solution actually exists for all times $t\in[0,T)$, so the maximal time horizon of existence is actually $T^L=T$, independently of the discretisation level $L$.

\noindent
{\bf Step 2.} energy estimates.\\
With $k_*$ defined by \eqref{eqn:mast} in case (i) and set to $k_*=0$ in case (ii), we use $v=\partial_t^{1+\alpha_{k_*}}u^L(s)$ as a test function in \eqref{eqn:Galerkin} (with $t$ replaced by $s$) and integrate for $s$ from $0$ to $t$. 
The resulting terms can be estimated as follows.

Inequality \eqref{eqn:coercivityI} 
\rerevision{with $\theta=1+\gamma_j-\alpha_{k_*}\in[0,1)$}
yields
\[
\begin{aligned}
&\int_0^t\langle d_j \partial_t^{2+\gamma_j} u^L(s), \partial_t^{1+\alpha_{k_*}}u^L(s)\rangle\, ds
= \int_0^t\langle d_j \partial_t^{2+\gamma_j} u^L(s), {_0I_t}^{1+\gamma_j-\alpha_{k_*}}\partial_t^{2+\gamma_j}u^L(s) \rangle\, ds\\
&\geq \cos ( \pi(1+\gamma_j-\alpha_{k_*})/2 ) \| \sqrt{d_j}\partial_t^{2+\gamma_j} u^L \|_{H^{-(1+\gamma_j-\alpha_{k_*})/2}(0,t;L^2_w(\Omega))}^2
\quad\mbox{ for }\gamma_j<\alpha_{k_*}
\end{aligned} 
\]
and \rerevision{with $\theta=1+\alpha_{k_*}-\alpha_k\in[0,1)$}
\[
\begin{aligned}
&\int_0^t\langle b_k \partial_t^{\alpha_k} \calAtil u^L(s), \partial_t^{1+\alpha_{k_*}}u^L(s)\rangle\, ds
=\int_0^t\langle b_k \partial_t^{\alpha_k} \calAtil^{1/2}u^L(s), \partial_t^{1+\alpha_{k_*}}\calAtil^{1/2}u^L(s)\rangle\, ds\\
&= \int_0^t\langle b_k {_0I_t}^{1+\alpha_{k_*}-\alpha_k}\partial_t^{1+\alpha_{k_*}} \calAtil^{1/2}u^L(s), \partial_t^{1+\alpha_{k_*}}\calAtil^{1/2}u^L(s) \rangle\, ds\\
&\geq \cos ( \pi(1+\alpha_{k_*}-\alpha_k)/2 ) \| \sqrt{b_k}\partial_t^{1+\alpha_{k_*}} \calAtil^{1/2}u^L \|_{H^{-(1+\alpha_{k_*}-\alpha_k)/2}(0,t;L^2_w(\Omega))}^2
\ \mbox{ for }k>{k_*}\,.
\end{aligned} 
\]
In particular for $k=K>{k_*}$ (which implies $1+\alpha_{k_*}-\alpha_K<1$; the case $k_*=K$ will be considered below), by \eqref{eqn:equivalence_neg} 
\rerevision{with $\theta=(1+\alpha_{k_*}-\alpha_K)/2\in[0,1)$}, noting that 
$\partial_t^{(1+\alpha_{k_*}+\alpha_K)/2} \calAtil^{1/2}u^L$ vanishes at $t=0$ and that $\partial_t^{(1+\alpha_{k_*}-\alpha_k)/2} \partial_t^{(1+\alpha_{k_*}+\alpha_K)/2}\calAtil^{1/2}u^L
=\partial_t^{1+\alpha_{k_*}} \calAtil^{1/2}u^L$,  
\begin{equation}\label{eqn:estbMupper}
\begin{aligned}
&\int_0^t\langle b_K \partial_t^{\alpha_K} \calAtil^{1/2}u^L(s), \partial_t^{1+\alpha_{k_*}}\calAtil^{1/2}u^L(s)\rangle\, ds\geq \rerevision{\underline{C}_{K,\alpha_{k_*}}} \| \partial_t^{(1+\alpha_{k_*}+\alpha_K)/2}\calAtil^{1/2}u^L \|_{L^2(0,t;L^2_w(\Omega))}^2
\end{aligned} 
\end{equation}
with 
\[
\rerevision{\underline{C}_{K,\alpha_{k_*}}}=\underline{C}((1+\alpha_{k_*}-\alpha_K)/2)\, b_K\, \cos ( \pi(1+\alpha_{k_*}-\alpha_K)/2 ). 
\] 

From 
\rerevision{the identity $\int_0^t\langle\partial_t v(s),v(s)\rangle ds =\frac12\|v(t)\|^2-\frac12\|v(0)\|^2$}
and Young's inequality in the case $\gamma_J=\alpha_{k_*}$ we get
\[
\begin{aligned}
&\int_0^t\langle d_J \partial_t^{2+\gamma_J} u^L(s), \partial_t^{1+\alpha_{k_*}}u^L(s)\rangle\, ds
=\int_0^t\langle d_J \Bigl(\partial_t\partial_t^{1+\alpha_{k_*}} u^L(s)
-\frac{s^{-\alpha_{k_*}}}{\Gamma(1-\alpha_{k_*})} u_2^L\Bigr), \partial_t^{1+\alpha_{k_*}}u^L(s)\rangle\, ds\\
&\geq\frac12\|\sqrt{d_J}\partial_t^{1+\alpha_{k_*}} u^L(t)\|_{L^2_w(\Omega)}^2
-\frac{t^{2(1-\alpha_{k_*})}\,\|\sqrt{d_J}u_2^L\|_{L^2_w(\Omega)}^2}{\Gamma(2-\alpha_{k_*})^2}
-\frac14 \|\sqrt{d_J}\partial_t^{1+\alpha_{k_*}}u^L\|_{L^\infty(0,t);L^2_w(\Omega)}^2\\
&\mbox{ if }\gamma_J=\alpha_{k_*}\in(0,1)
\end{aligned} 
\]
\[
\begin{aligned}
&\int_0^t\langle d_J \partial_t^{2+\gamma_J} u^L(s), \partial_t^{1+\alpha_{k_*}}u^L(s)\rangle\, ds
=\frac12\|\sqrt{d_J}\partial_t^{1+\alpha_{k_*}} u^L)(t)\|_{L^2_w(\Omega)}^2
-\frac12\|\sqrt{d_J} u^L_{1+\alpha_{k_*}}\|_{L^2_w(\Omega)}^2\\
&\mbox{ if }\gamma_J=\alpha_{k_*}\in\{0,1\}
\end{aligned} 
\]
and 
\[
\begin{aligned}
&\int_0^t\langle b_{k_*} \partial_t^{\alpha_{k_*}} \calAtil^{1/2}u^L(s), \partial_t^{1+\alpha_{k_*}}\calAtil^{1/2}u^L(s)\rangle\, ds\\
&=\int_0^t\langle b_{k_*} \partial_t\partial_t^{\alpha_{k_*}} \calAtil^{1/2}u^L(s)
-\frac{s^{-\alpha_{k_*}}}{\Gamma(1-\alpha_{k_*})} \calAtil^{1/2}u_1^L, \partial_t^{\alpha_{k_*}}\calAtil^{1/2}u^L(s)\rangle\, ds\\
&\geq\frac12\|\sqrt{b_{k_*}}\partial_t^{\alpha_{k_*}} \calAtil^{1/2}u^L)(t)\|_{L^2_w(\Omega)}^2
-\frac{t^{2(1-\alpha_{k_*})}\,\|\sqrt{b_{k_*}}\calAtil^{1/2}u_1^L\|_{L^2_w(\Omega)}^2}{\Gamma(2-\alpha_{k_*})^2} 
\\&\qquad
-\frac14 \|\sqrt{b_{k_*}}\partial_t^{\alpha_{k_*}}\calAtil^{1/2}u^L\|_{L^\infty(0,t);L^2_w(\Omega)}^2
\ \mbox{ if }\alpha_{k_*}\in(0,1)
\,,
\end{aligned} 
\]
\[
\begin{aligned}
&\int_0^t\langle b_{k_*} \partial_t^{\alpha_{k_*}} \calAtil^{1/2}u^L(s), \partial_t^{1+\alpha_{k_*}}\calAtil^{1/2}u^L(s)\rangle\, ds\\
&=\frac12\|\sqrt{b_{k_*}}\partial_t^{\alpha_{k_*}} \calAtil^{1/2}u^L)(t)\|_{L^2_w(\Omega)}^2
-\frac12\|\sqrt{b_{k_*}} \calAtil^{1/2}u^L_{\alpha_{k_*}}\|_{L^2_w(\Omega)}^2
\ \mbox{ if }\alpha_{k_*}\in\{0,1\}
\,.
\end{aligned} 
\]

Hence in the case $k_*=K$, where due to our assumption that $\alpha_K>0$ at the beginning of the proof, $\partial_t^{\alpha_K} \calAtil^{1/2}u^L(0)=0$,
\[
\begin{aligned}
&\sup_{t'\in(0,t)}\int_0^{t'}\langle b_K \partial_t^{\alpha_K} \calAtil^{1/2}u^L(s), \partial_t^{1+\alpha_{k_*}}\calAtil^{1/2}u^L(s)\rangle\, ds\\
&\geq \frac{b_K}{4}\|\partial_t^{\alpha_K}\calAtil^{1/2}u^L\|_{L^\infty(0,t);L^2_w(\Omega)}^2
-\frac{t^{2(1-\alpha_{K})}\,\|\sqrt{b_{K}}\calAtil^{1/2}u_1^L\|_{L^2_w(\Omega)}^2}{\Gamma(2-\alpha_{K})^2}.
\end{aligned} 
\]

Since the difference $1+\alpha_{k_*}-\alpha_k$ is larger than one for $k< k_*$, 
we have to employ \eqref{eqn:boundednessI} in that case \rerevision{(for this purpose, note that still, according to our assumptions, $\theta=(1+\alpha_{k_*}-\alpha_k)/2\in[0,1]$)}.
From this and \eqref{eqn:equivalence_neg}, 
\rerevision{again with $\theta=(1+\alpha_{k_*}-\alpha_k)/2\in[0,1]$,}
we obtain
\begin{equation}\label{eqn:estbmlower}
\begin{aligned}
&\int_0^t\langle b_k \partial_t^{\alpha_k} u^L(s), \partial_t^{1+\alpha_{k_*}}u^L(s)\rangle\, ds\\
&= \int_0^t\langle b_k \, {_0I_t}^{1+\alpha_{k_*}-\alpha_k}\partial_t^{1+\alpha_{k_*}} \calAtil^{1/2}u^L(s), \partial_t^{1+\alpha_{k_*}}\calAtil^{1/2}u^L(s) \rangle\, ds\\
&\geq -C(1+\alpha_{k_*}-\alpha_k) \|\sqrt{b_k}\partial_t^{1+\alpha_{k_*}} \calAtil^{1/2}u^L\|_{H^{-(1+\alpha_{k_*}-\alpha_k)/2}(0,t;L^2_w(\Omega))}^2\\
&\geq -C(1+\alpha_{k_*}-\alpha_k)\overline{C}((1+\alpha_{k_*}-\alpha_k)/2)\|\sqrt{b_k}\partial_t^{(1+\alpha_{k_*}+\alpha_k)/2}\calAtil^{1/2}u^L\|_{L^2(0,t;L^2_w(\Omega))}^2\\
&\mbox{ for }k<{k_*}\,.
\end{aligned} 
\end{equation}

\rerevision{
Note that in case $k_*<K$, the potentially negative contributions arising from the $k<k_*$ terms \eqref{eqn:estbmlower} are dominated by means of the leading energy term from \eqref{eqn:Eb}, that is, $\| \partial_t^{(1+\alpha_{k_*}+\alpha_K)/2}\calAtil^{1/2}u^L \|_{L^2(0,t;L^2_w(\Omega))}^2$ from \eqref{eqn:estbMupper}, due to the fact that $(1+\alpha_{k_*}+\alpha_k)/2< (1+\alpha_{k_*}+\alpha_K)/2$ for $k<{k_*}$.
}

\rerevision{To achieve this also in case $k_*=K$ with $\alpha_{K-1}<\gamma_J<\alpha_K$, we additionally test with $v=\partial_t^{1+\alpha_K-\epsilon}u^L(s)$ in \eqref{eqn:Galerkin} where $0<\epsilon<\alpha_K-\max\{\alpha_{K-1},\gamma_J\}$.
Repeating the above estimates with $\alpha_{k_*}$ replaced by $\alpha_K-\epsilon$ one sees that 
the leading energy term then is $\| \partial_t^{\alpha_K+(1-\epsilon)/2}\calAtil^{1/2}u^L \|_{L^2(0,t;L^2_w(\Omega))}^2$ from \eqref{eqn:estbMupper}. This still dominates the potentially negative contributions arising from the $k<K$ terms according to \eqref{eqn:estbmlower} (with $\alpha_K-\epsilon$ in place of $\alpha_{k_*}$), since $(1+\alpha_K-\epsilon+\alpha_k)/2< \alpha_K+(1-\epsilon)/2$ for $k<K$.\\
Note that in the case $\gamma_J=\alpha_K=\alpha_0$, no potentially negative terms \eqref{eqn:estbmlower} arise.}

Finally, the right hand side term can be estimated by means of \eqref{eqn:equivalence_neg} 
\rerevision{with $\theta=\alpha_{k_*}-\gamma_J\in[0,1)$}
and Young's inequality  
\[
\begin{aligned}
&\int_0^t \langle r(s), \partial_t^{1+\alpha_{k_*}}u^L(s)\rangle\, ds
\leq\|r\|_{H^{\alpha_{k_*}-\gamma_J}(0,T,L^2_w(\Omega))}
\|\partial_t^{1+\alpha_{k_*}}u^L(s)\|_{H^{-(\alpha_{k_*}-\gamma_J)}(0,t,L^2_w(\Omega))}\\
&\leq \overline{d} \|\partial_t^{1+\gamma_J}u^L(s)\|_{L^2(0,t,L^2_w(\Omega))}^2
+\frac{1}{4\overline{d}\overline{C}(\alpha_{k_*}-\gamma_J)^2} \|r\|_{H^{\alpha_{k_*}-\gamma_J}(0,T,L^2_w(\Omega))}^2\\
&\rerevision{\leq \int_0^t \underline{\mathcal{E}}_\gamma[u^L](s)\, ds + C \|r\|_{H^{\alpha_{k_*}-\gamma_J}(0,T,L^2_w(\Omega))}^2}
\end{aligned} 
\]
\rerevision{with some constants $\overline{d}, C>0$.}

Altogether we arrive at the following estimate 
\begin{equation}\label{eqn:enest0}
\mathcal{E}_\gamma[u^L](t)+\mathcal{E}_\alpha[u^L](t)\leq 
\rerevision{\int_0^t \underline{\mathcal{E}}_\gamma[u^L](s)\, ds} 
+\mbox{rhs}[u^L](t) + 
\mbox{rhs}^L_{0r}(t),
\end{equation}
where $\mathcal{E}_\gamma$, \rerevision{$\underline{\mathcal{E}}_\gamma$}, $\mathcal{E}_\alpha$ are defined as in \eqref{eqn:Ea}, \rerevision{\eqref{eqn:ulEa}}, \eqref{eqn:Eb},
\[
\begin{aligned}
&\mbox{rhs}[u](t)=
\sum_{k=0}^{k_*-1}C(1+\alpha_{k_*}-\alpha_k)\overline{C}((1+\alpha_{k_*}-\alpha_k)/2)\|\sqrt{b_k}\partial_t^{(1+\alpha_{k_*}+\alpha_k)/2}\calAtil^{1/2}u\|_{L^2(0,t;L^2_w(\Omega))}^2\\
&\mbox{rhs}^L_{0,r}(t)=
1_{\gamma_J=\alpha_{k_*}\in(0,1)} \frac{t^{2(1-\alpha_{k_*})}\,\|\sqrt{d_J}u_2^L\|_{L^2_w(\Omega)}^2}{\Gamma(2-\alpha_{k_*})^2}
+1_{\alpha_{k_*}\in(0,1)} \frac{t^{2(1-\alpha_{k_*})}\,\|\sqrt{b_{k_*}}\calAtil^{1/2}u_1^L\|_{L^2_w(\Omega)}^2}{\Gamma(2-\alpha_{k_*})^2}\\
&\hspace*{2cm}+\frac{1}{2\tilde{d}\overline{C}(\alpha_{k_*}-\gamma_J)^2} \|r\|_{H^{\alpha_{k_*}-\gamma_J}(0,T,L^2_w(\Omega))}^2
\end{aligned}
\]
where the boolean variable $1_B$ is equal to one if $B$ is true and vanishes otherwise.

We now substitute $\mathcal{E}_\gamma[u](t)$, $\mathcal{E}_\alpha[u](t)$ by lower bounds containing only an estimate from below of their leading order terms  $\underline{\mathcal{E}}_\gamma[u](t)$ as defined in \eqref{eqn:ulEa} and 
\begin{equation}\label{eqn:ulEb}
\underline{\mathcal{E}}_\alpha[u](t) :=
\begin{cases}
\rerevision{\underline{C}_{K,\alpha_{k_*}}} \| \partial_t^{(1+\alpha_{k_*}+\alpha_K)/2}\calAtil^{1/2}u \|_{L^2(0,t;L^2_w(\Omega))}^2
&\mbox{ if }k_*<K\\
\rerevision{\underline{C}_{K,\alpha_K-\epsilon} \| \partial_t^{\alpha_K+(1-\epsilon)/2}\calAtil^{1/2}u \|_{L^2(0,t;L^2_w(\Omega))}^2}
&\mbox{ if }k_*=K.
\end{cases}
\end{equation}
We then obtain from \eqref{eqn:enest0}
\begin{equation}\label{eqn:enest1}
\underline{\mathcal{E}}_\gamma[u^L](t)+\underline{\mathcal{E}}_\alpha[u^L](t)\leq 
\rerevision{\int_0^t \underline{\mathcal{E}}_\gamma[u^L](s)\, ds} 
+\mbox{rhs}[u^L](t) + \mbox{rhs}_{0,r}^L(t)
+ \widetilde{\mbox{rhs}}_0^L(t)
\end{equation}
where 
\[
\widetilde{\mbox{rhs}}_0^L(t)=
1_{k_*=K,\,\alpha_K\in(0,1)}\frac{t^{2(1-\alpha_{K})}\,\|\sqrt{b_{K}}\calAtil^{1/2}u_1^L\|_{L^2_w(\Omega)}^2}{\Gamma(2-\alpha_{K})^2}\,.
\]
\rerevision{
It is readily checked that 
$
\mbox{rhs}[u](t)\leq c_{b,\alpha}(t)\underline{\mathcal{E}}_\alpha[u](t)
$
for 
\[
c_{b,\alpha}(t)= \frac{1}{\underline{C}_{K,\alpha_{k_*}}}\sum_{k=0}^{k_*-1} C(1+\alpha_{k_*}-\alpha_k)\overline{C}((1+\alpha_{k_*}-\alpha_k)/2) |b_k| \|_0I_t^{(\alpha_K-\alpha_k)/2}\|_{L^2(0,t)\to L^2(0,t)} 
.\]
Therefore, under the smallness assumption \eqref{smallnessbk}, that is, $c_{b,\alpha}(T)<1$ 
estimate \eqref{eqn:enest1} yields
\[
\begin{aligned}
&\underline{\mathcal{E}_\gamma}[u^L](t)+(1-c_{b,\alpha}(T))\underline{\mathcal{E}}_\alpha[u^L](t)
\leq \tilde{C} \bigl(\rerevision{\int_0^t \underline{\mathcal{E}}_\gamma[u^L](s)\, ds} 
+\|u_0\|_{\dot{H}^1(\Omega)}^2 + C_0
+\|r\|_{H^{\alpha_{k_*}-\gamma_J}(0,T,L^2_w(\Omega))}^2\bigr)
\end{aligned}
\]
for some $\tilde{C}$ independent of $t$, $L$, with $C_0$ as in \eqref{eqn:C0}.
Applying Gronwall's inequality and }
re-inserting into \eqref{eqn:enest0} we end up with 
\begin{equation}\label{eqn:enest_sumfrac_L}
\begin{aligned}
&\mathcal{E}_\gamma[u^L](t)+\mathcal{E}_\alpha[u^L](t)
\leq C(T) \bigl(\|u_0\|_{\dot{H}^1(\Omega)}^2 + C_0 +\|r\|_{H^{\alpha_{k_*}-\gamma_J}(0,T,L^2_w(\Omega))}^2\bigr).
\end{aligned}
\end{equation}

\noindent
{\bf Step 3.} weak limits.\\
As a consequence of \eqref{eqn:enest_sumfrac_L}, the Galerkin solutions $u^L$ are uniformly bounded in $U_\gamma\cap U_\alpha$. Therefore the sequence $(u^L)_{L\in\mathbb{N}}$ has a weakly(*) convergent subsequence $(u^{L_k})_{k\in\mathbb{N}}$ with limit $u\in U_\gamma\cap U_\alpha$.
To see that $u$ satisfies \eqref{eqn:sumfrac_wave_acou_ibvp_timeint}, we revisit \eqref{eqn:Galerkin}, integrate it with respect to time and take the $L^2(0,T)$ product with arbitrary smooth test functions $\psi$ to conclude that  
\[
\begin{aligned}
&\int_0^T\psi(t)\Bigl\{\langle \sum_{j=0}^J d_j \Bigl(\partial_t^{1+\gamma_j} u(t)-\frac{t^{1-\gamma_j}}{\Gamma(2-\gamma_j)}u_2\Bigr),v\rangle_{L^2_w(\Omega)} \\
&\qquad\qquad+ \langle\sum_{k=0}^K b_k \Bigl({_0I_t}^{1-\alpha_k} \calAtil u(t) - \frac{t^{1-\alpha_k}}{\Gamma(2-\alpha_k)}\calAtil u_0\Bigr) - r(t),v\rangle_{\dot{H}(\Omega)^*,\dot{H}(\Omega)}\Bigr\}\, dt \\
&= \int_0^T\psi(t)\Bigl\{\langle \sum_{j=0}^J d_j \partial_t^{1+\gamma_j} (u-u^{L_k})(t),v\rangle_{L^2_w(\Omega)}\\ 
&\qquad\qquad + \langle\sum_{k=0}^K b_k {_0I_t}^{1-\alpha_k} \calAtil (u-u^{L_k})(t),v\rangle_{\dot{H}(\Omega)^*,\dot{H}(\Omega)}\Bigr\}\, dt  
\end{aligned}
\]
for all $v\in \mbox{span}(\varphi_1,\ldots, \varphi_K)$, $K\leq L_k$ and any $\psi\in C_c^\infty(0,T)$. Taking the limit $k\to\infty$ by using the previously mentioned weak(*) limit of $u^{L_k}$ 
we conclude 
\[
\begin{aligned}
&\int_0^T\psi(t)\Bigl\{\langle \sum_{j=0}^J d_j \Bigl(\partial_t^{1+\gamma_j} u(t)-\frac{t^{1-\gamma_j}}{\Gamma(2-\gamma_j)}u_2\Bigr),v\rangle_{L^2_w(\Omega)} \\
&\qquad\qquad+ \langle\sum_{k=0}^K b_k \Bigl({_0I_t}^{1-\alpha_k} \calAtil u(t) - \frac{t^{1-\alpha_k}}{\Gamma(2-\alpha_k)}\calAtil u_0\Bigr) - r(t),v\rangle_{\dot{H}(\Omega)^*,\dot{H}(\Omega)}\Bigr\}\, dt 
=0
\end{aligned}
\]
for all $v\in \mbox{span}(\varphi_1,\ldots, \varphi_K)$ with arbitrary $K\in\mathbb{N}$, and any $\psi\in C_c^\infty(0,T)$. 
Therefore the weak limit $u$ indeed satisfies the time integrated {\sc pde} and we have proven the existence part of the theorem. 

To verify the initial conditions, we use the fact that due to our assumption $\alpha_K>0$, $U_\alpha$ continuously embeds into $C([0,T];\dot{H}(\Omega))$ and therefore $u(0)=u_0$ is attained in an $\dot{H}(\Omega)$ sense. 
Also, $U_\gamma$ continuously embeds into $C^1([0,T];L^2_w(\Omega))$ in case 
$\gamma_J<\alpha_{k_*}$ or $\gamma_J>0$. 
In the remaining case $\gamma_J=\alpha_{k_*}=0$, attainment of $u_t(0)=u_1$ in an $L^2_w(\Omega)$ sense can be shown analogously to the conventional second order wave equation, cf., e.g., \cite[Theorem 3 in Section 7.2]{Evans:2010}. 

Moreover, taking weak limits in the energy estimate \eqref{eqn:enest_sumfrac_L}
together with weak lower semicontinuity of the norms contained in the definitions
of $\mathcal{E}_\gamma$, $\mathcal{E}_\alpha$, implies \eqref{eqn:enest_sumfrac}.

\noindent
{\bf Step 4.} uniqueness.\\
The manipulations carried out in Step 2. of the proof are also feasible with the Galerkin approximation $u^L$ replaced by a solution $u$ of the {\sc pde} itself, and lead to the energy estimate \eqref{eqn:enest_sumfrac} independently of the Galerkin approximation procedure. 
From this we conclude that the  {\sc pde} with vanishing right hand side and initial data only has the zero solution, which due to linearity of the {\sc pde} implies uniqueness.
\end{proof}

\rerevision{
\begin{corollary}\label{cor:global}
Under the assumptions of Theorem~\ref{thm:sumfrac} 
 with $b_k=0$, $k=0,\ldots,k_*-1$, the assertion extends to hold globally in time, that is, for all $t\in(0,\infty)$. If additionally, $r=0$, then the energy is globally bounded 
\[
\mathcal{E}_\gamma[u](t)+\mathcal{E}_\alpha[u](t) \leq 
C \bigl(\|u_0\|_{\dot{H}^1(\Omega)}^2 + C_0\bigr)\,, \quad t\in(0,\infty)
\]
with some constant $C>0$ independent of time. 
\end{corollary}
}
\begin{remark}
Using the so-called multinomial (or multivariate) Mittag-Leffler functions and separation of variables in principle enables a solution representation by separation of variables, cf. \cite[Theorem 4.1]{LuchkoGorenflo:1999}.
However, proving convergence of the infinite sums in this representation would require extensive estimates of these Mittag-Leffler functions. Thus, in order to keep the exposition transparent without having to introduce too much additional machinery, we chose to remain with the energy argument in the proof above. 
\end{remark}

In order to establish sufficient regularity of the highest order term $\partial_t^{2+\gamma_J} u^L$ so that also the original equation \eqref{eqn:sumfrac_wave_acou_ibvp} holds, we make use of the {\sc pde} itself. 
\begin{corollary}
Under the conditions of Theorem~\ref{thm:sumfrac} 
the solution $u$ of \eqref{eqn:sumfrac_wave_acou_ibvp_timeint} 
satisfies $\partial_t^{2+\gamma_J} u\in L^2(0,T;\dot{H}(\Omega)^*)$ and is therefore also a solution to the original {\sc pde} \eqref{eqn:sumfrac_wave_acou_ibvp} in an $L^2(0,T;\dot{H}(\Omega)^*)$ sense.
\end{corollary}
\begin{proof}
From the Galerkin approximation \eqref{eqn:Galerkin} of \eqref{eqn:sumfrac_wave_acou_ibvp}, due to invariance of the space $\mbox{span}(\varphi_1,\ldots, \varphi_L)$ we can substitute $v$ by $\calAtil^{-1/2}v$ there and move $\calAtil^{-1/2}$ to the left hand side of the inner product by taking the adjoint
to arrive at 
\begin{equation}\label{eqn:Galerkin_lower}
\begin{aligned}
&\langle \sum_{j=0}^J \partial_t^{2+\gamma_j}\calAtil^{-1/2}[d_j  u^L(t)] + \sum_{k=0}^K \partial_t^{\alpha_k} \calAtil^{-1/2}[b_k  
\calAtil u^L(t)] - \calAtil^{-1/2}r(t),v\rangle_{L^2_w(\Omega)} = 0 \\
&\quad t\in(0,T) \quad \mbox{ for all }v\in\mbox{span}(\varphi_1,\ldots, \varphi_L)\,.
\end{aligned}
\end{equation}
From Theorem~\ref{thm:sumfrac} we know that 
$\tilde{r}=-\sum_{k=0}^K \partial_t^{\alpha_k} \calAtil^{-1/2}[b_k  \calAtil u(t)] + \calAtil^{-1/2}r(t)$ is contained in $L^2(0,T;L^2_w(\Omega))$. 
Thus \eqref{eqn:Galerkin_lower} can be viewed as a Galerkin discretisation of a Volterra integral equation for $\zeta =  \partial_t^{2+\gamma_J}\calAtil^{-1/2} u(t)$ with $\tilde{r}$ as right hand side. The corresponding coefficient vector $\underline{\zeta}^L$ of the solution $\zeta^L$ is therefore bounded by 
\[
\|\underline{\zeta}^L\|_{\ell^2(\mathbb{R}^L)}
\leq C(\|\tilde{r}\|_{L^2(0,T;L^2_w(\Omega))}+\|\zeta(0)\|_{L^2_w(\Omega)})
\]
with $C$ independent of $L$. From this we deduce 
\[ \|\zeta^L\|_{L^2(0,T;L^2_w(\Omega))}
\leq C(\|\tilde{r}\|_{L^2(0,T;L^2_w(\Omega))}+\|\zeta(0)\|_{L^2_w(\Omega)})
\]
and thus, by taking the limit $L\to\infty$ the same bound for $\|\zeta^L\|_{L^2(0,T;L^2_w(\Omega))}$ itself.
\end{proof}

%% file: numerics_forward.tex
\subsection{Numerical solution of the forward problem}

There are several methods available for solving multi-term
fractional ordinary differential equations.
For our situation we require a solver that is accurate over the entire
time interval $[0,\infty)$ since one problem we consider
is the recovery of the components of the operator from large time values only.
Our method of choice here is to convert the system with multiple
fractional derivatives to a single 
{\sc ode}
of the form  $D_t^\gamma u = A u$
where $A$ is a square matrix.
The seminal starting point here is the paper by Diethelm and Ford,
\cite{DiethelmFord:2004} which allows one to consider the equation
$D^{\alpha_0} y(t) = f\bigl(t,y(t),D^{\alpha_1},\;\ldots\;,D^{\alpha_n}\bigr)$
subject to $y^{(k)}(0) = y_0^k$ for $k=0,\,1,\;\ldots\;\lceil\alpha_0\rceil-1$,
and its conversion to the form $D^\gamma Y(t) = g\bigl(t,Y(t)\bigr)$,
with $Y(0) = Y_0$.
Existence and uniquness results are obtained and a numerical method for
the solution formulated.
This was later refined by Garrappa and Popolizio,
\cite{GarrappaPopolizio:2018} in the case of linear systems
into a solver whose accuracy approaches machine precision.
This is based on using the Mittag-Leffler function with a matrix argument $A$
that uses different approaches depending of the location of the
eigenvalues of $A$ in the complex plane and combined into an algorithm that
works for matrices with any eigenvalue locations.
See also \cite{Popolizio:2018}.
We describe this above approach briefly since it comes with a certain
restriction that should be clarified.

The following result is proven in \cite{Popolizio:2018}:
\begin{theorem}
Consider the initial value problem
\begin{equation}\label{eqn:ode_ivp}
D^{\alpha_0} = f\bigl(t,y(t),D^{\alpha_1}y(t),\,\ldots\,
D^{\alpha_{n-1}}y(t)\bigr),\qquad
y^{(j)}(0) = y_0{(j)}, j = 0,\;\ldots\,,\,[\alpha_0]-1
\end{equation}
where $\alpha_0> \alpha_1 >\, \ldots\, \alpha_{k-1}$ and $0<\alpha_{k-1}<1$.
Assume that each $\alpha_j\in \mathbb{Q}$ and let $M$ 
be the least
common denominator of $\alpha_0,\,\ldots,\alpha_{k-1}$.
Set $\gamma=1/M$ and $N=M\alpha_0$.
Then the initial value problem \eqref{eqn:ode_system} is equivalent to
the system of equations
\begin{equation}\label{eqn:ode_system}
D^\gamma y_0(t) = y_1(t),\ D^\gamma y_1(t) = y_2(t), \ 
	D^\gamma y_{N-1}(t) = f(t,y_0(t),y_{\alpha_1/\gamma}(t),\ldots
	\,y_{\alpha_{k-1}/\gamma}(t)
\end{equation}
together with the initial conditions
\begin{equation}\label{eqn:ode_system_init_conditions}
	y_j(0) = y_0^{(j/M)} \quad \mbox{if } j/M \in \mathbb{Z}_0,
	\quad \mbox{else }\ 0.
\end{equation}
\end{theorem}
As a result of the above theorem the linear system can be reformulated to
the form
\begin{equation}\label{eqn:ode_matrix_system}
D^\gamma Y(t) = A Y(t) + 
e_N\,f(t),\qquad Y(0)=Y_0,
\quad e_N = (0,0,\ldots,1)^T\in \mathbb{R}^N.
\end{equation}
where the coefficient matrix $A\in \mathbb{R}^{N\times N}$ is of 
companion type.

Then it is possible to express the solution to \eqref{eqn:ode_matrix_system}
in terms of Mittag-Leffler functions with matrix arguments as
\begin{equation}\label{eqn:ode_matrix_ML}
Y(t) = E_{\alpha,1}(t^\gamma A) Y_0 
+ \int_0^t (t-s)^{\gamma-1} E_{\gamma,\gamma}\bigl((t-s)^\gamma A\bigr)
	e_N\,f(s)\,ds.
\end{equation}
For general matrix functions the Schur-Parlett algorithm is the method of
choice and to this end the final component needed is an accurate evaluation
of the matrix Mittag-Leffler function and its derivatives.
This is both nontrivial and quite technical and we simply refer to the above
references, \cite{GarrappaPopolizio:2018,Popolizio:2018},
and also to \cite{Garrappa:2015} for the details.
Code to implement what is needed for this step is available from
Roberto Garrappa's homepage \cite{Garrappa}.

Note that in our case we have $\alpha_0=2$ and also the requirement that
each of the exponents $\alpha_j = q_j/p_j$ must be rational numbers.
The above restrictions can lead to a quite large system of size
$\mbox{lcm}(p_1,\,\ldots,p_{k-1})$ and hence a small value of $\gamma$.
However, we use this solver only for the purpose of providing simulated data
and our inversion techniques are independent of its use.

\revision{
Since we are using asymptotic formulas it is essential to use a reliable and accurate solver that is capable to cover a wide range of time values. The code from \cite{Garrappa:2015,GarrappaPopolizio:2018,Popolizio:2018} provably fulfills these requirements and our numerical experiments confirmed this.
}

%% file: resolvent.tex
\subsection{Resolvent equations}

Separation of variables based on the eigensystem of $\mathcal{A}$ yields a solution representation
\[
u(x,t)=\sum_{n=1}^\infty \Bigl(
\sigma(t)w_{fn}(t)\langle f,\varphi_n\rangle 
+w_{2n}(t)\langle u_2,\varphi_n\rangle 
+w_{1n}(t)\langle u_1,\varphi_n\rangle 
+w_{0n}(t)\langle u_0,\varphi_n\rangle
\Bigr)\varphi_n(x)
\]
and convergence of the sums is guaranteed by the energy estimates from Theorem~\ref{thm:sumfrac}.
Here $w_{fn}$, $w_{2n}$, $w_{1n}$, $w_{0n}$ denote the solutions of the
resolvent equation
\[
\begin{aligned}
&w_{n,tt}+c^2\lambda_n w_n
+\sum_{j=1}^J d_j D_t^{2+\gamma_j}w_n
+\sum_{k=1}^N b_k\lambda_n^{\beta_k} D_t^{\alpha_k} w_n = 
\begin{cases}1\ \mbox{ for }w_{fn}\\0\ \mbox{ else}\end{cases}\\
&
w_{n,tt}(0)=\begin{cases}1\ \mbox{ for }w_{2n}\\0\ \mbox{ else}\end{cases}, \quad
w_{n,t}(0)=\begin{cases}1\ \mbox{ for }w_{1n}\\0\ \mbox{ else}\end{cases}, \quad
w_{n}(0)=\begin{cases}1\ \mbox{ for }w_{0n}\\0\ \mbox{ else}\end{cases}.
\end{aligned}
\]

From the Laplace transformed resolvent equations we obtain the resolvent solutions
\[ \begin{aligned}
&\hat{w}_{fn}(s)=\frac{1}{
\omega(s;\lambda_n)}
\,, \quad
\hat{w}_{0n}(s)=
\frac{s+\sum_{j=1}^J d_j s^{\gamma_j+1}+\sum_{k=1}^N s^{\alpha_k-1} b_{k}\lambda_n^{\beta_{k}}}{
\omega(s;\lambda_n)}
=\frac{\omega(s,\lambda_n)-c^2\lambda_n}{s\omega(s;\lambda_n)}
\,, \\
&\hat{w}_{1n}(s)=\frac{1+\sum_{j=1}^J d_j s^{\gamma_j}}{
\omega(s;\lambda_n)}
\,, \quad
\hat{w}_{2n}(s)=\frac{\sum_{j=1}^J d_j s^{\gamma_j-1}}{
\omega(s;\lambda_n)}
\end{aligned}\]
with 
\[
\omega(s;\lambda_n)=s^2+c^2\lambda_n+\sum_{j=1}^J d_j s^{\gamma_j+2}+\sum_{k=1}^N s^{\alpha_k} b_{k}\lambda_n^{\beta_{k}},
\]
where $\hat{w}$ denotes the Laplace transform of $w$.

Formally we can write
\[
\hat{u}(x,s)=\sum_{n=1}^\infty \Bigl(
\hat{\sigma}(s)\hat{w}_{fn}(s)\langle f,\varphi_n\rangle 
+\hat{w}_{2n}(s)\langle u_2,\varphi_n\rangle 
+\hat{w}_{1n}(s)\langle u_1,\varphi_n\rangle 
+\hat{w}_{0n}(s)\langle u_0,\varphi_n\rangle
\Bigr)\varphi_n(x)
\]
and therefore
\[
\widehat{B_iu}(s)=\sum_{n=1}^\infty \Bigl(
\hat{\sigma}(s)\hat{w}_{fn}(s)\langle f,\varphi_n\rangle 
+\hat{w}_{2n}(s)\langle u_2,\varphi_n\rangle 
+\hat{w}_{1n}(s)\langle u_1,\varphi_n\rangle 
+\hat{w}_{0n}(s)\langle u_0,\varphi_n\rangle
\Bigr)B_i\varphi_n\,.
\]

%% file: uniqueness.tex
\section{Uniqueness approaches}\label{sec:uniqueness}
Choose, for simplicity the excitations $u_{0i}=0$, $u_{1i}=0$, $u_{2i}=0$, $f_i\neq0$ and assume that 
for two models given by 
\begin{equation}\label{eqn:model}
\begin{aligned}
&N, \ J, \ c, \ \mathcal{A}, \ b_{k}, \ \alpha_k, \ \beta_{k}, \ d_j, \ \gamma_j,\ B_i,\\ 
&\tilde N, \ \tilde J, \ \tilde c, \ \tilde{\mathcal{A}}, \ \tilde{b}_{k}, \ \tilde{\alpha}_k, \ \tilde \beta_{k}, \ \tilde d_j, \ \tilde \gamma_j,\ \tilde B_i,
\end{aligned}
\end{equation}
we have equality of the observations
\[
B_i u(t) = \tilde B_i \tilde u(t)\, \quad t\in(0,T)
\]
from two possibly different space parts of the excitation $f_i$, $\tilde{f}_i$, and possibly using two different test specimen characterised by $\mathcal{A}$, $\tilde{\mathcal{A}}$,
while the temporal part $\sigma$ of the excitation is assumed to be the same in both experiments, and its Laplace transform to vanish nowhere $\hat{\sigma}(s)\neq0$ for all $s$.
By analyticity (due to the multinomial Mittag-Leffler functions reresentation cf. \cite[Theorem 4.1]{LuchkoGorenflo:1999})
we get equality for all $t>0$, thus equality of the Laplace transforms
\begin{equation}\label{eqn:BaBb}
\begin{aligned}
&\sum_{n=1}^\infty \frac{\langle f_{i},\varphi_n\rangle B_i\varphi_n}{
s^2+c^2\lambda_n+\sum_{j=1}^J d_j s^{\gamma_j+2}+\sum_{k=1}^N s^{\alpha_k} b_{k}(\lambda_n)^{\beta_{k}}}\\
&\qquad=
\sum_{n=1}^\infty \frac{\langle {\tilde f}_{i},{\tilde \varphi}_n\rangle {\tilde B}_i {\tilde \varphi}_n}{
s^2+{\tilde c}^2 {\tilde \lambda}_n+\sum_{j=1}^{\tilde{J}} \tilde{d}_j s^{\tilde{\gamma}_j+2}+\sum_{k=1}^{\tilde N} s^{{\tilde \alpha}_k} {\tilde b}_{k}({\tilde \lambda}_n)^{{\tilde \beta}_{k}}}, \quad i=1,\ldots I\,, \quad s\in\mathbb{C}_\theta
\end{aligned}
\end{equation}
in some sector $\mathbb{C}_\theta$ of the complex plane,
due to our assumption that $\hat{\sigma}(s)=\hat{\tilde{\sigma}}(s)$ and $\hat{\sigma}(s)\neq0$ for all $s$.
We have to disentangle the sum over $n$ with the asymptotics for $s\to\infty$ in order to recover the unknown quantities in \eqref{eqn:model} or at least some of them.

\subsection{Smooth excitations}

We make use of the Weyl estimate $\lambda_n \approx C_d\, n^{2/d}$ in
$\mathbb{R}^d$ as well as additional regularity of the excitation to approximate
the sum over $n$ \eqref{eqn:BaBb} of fractions by a single fraction of the form
\[
\sum_{n=1}^\infty \frac{\langle f_{i},\varphi_n\rangle B_i\varphi_n}{\omega(s,\lambda_n)}
\approx\frac{\sum_{n=1}^\infty\langle f_{i},\varphi_n\rangle B_i\varphi_n}{\omega(s,\lambda_1)}
\mbox{ with } \omega(s,\lambda) = s^2+c^2\lambda+\sum_{j=1}^J d_j s^{\gamma_j+2}+\sum_{k=1}^N s^{\alpha_k} b_{k}\lambda^{\beta_{k}}
\] 
for large real positive $s$.
We focus exposition on the source excitation case $u_\ell=0$, $\ell=0,1,2$, $f\neq0$ from \eqref{eqn:BaBb}, the other cases follow similar. Also, since we only need one excitation here $I=1$, we just skip the subscript $i$.
Indeed, for $s\in\mathbb{R}$, $s\geq1$ we can estimate the difference
\[
\begin{aligned}
D:=&\left|
\frac{\sum_{n=1}^\infty\langle f,\varphi_n\rangle B\varphi_n}{\omega(s,\lambda_1)}-
\sum_{n=1}^\infty \frac{\langle f,\varphi_n\rangle B\varphi_n}{\omega(s,\lambda_n)}
\right|\\
&=\left|\sum_{n=1}^\infty
\frac{\langle f,\varphi_n\rangle B\varphi_n\,(c^2(\lambda_n-\lambda_1)+\sum_{k=1}^N s^{\alpha_k} b_{k}((\lambda_n)^{\beta_{k}}-(\lambda_1)^{\beta_{k}}))}{\omega(s,\lambda_n)\, \omega(s,\lambda_1)}\right|\\
&\leq \sup_{n\in\mathbb{N}} \{\lambda_n^{-\nu}|B\varphi_n|\}\, 
\|f\|_{\dot{H}^{2\mu+2\nu}(\Omega)}
\left(\sum_{n=1}^\infty\left(
\frac{c^2\lambda_n+\sum_{k=1}^N s^{\alpha_k} b_{k}(\lambda_n)^{\beta_{k}}}{
\lambda_n^\mu s^2(s^2+c^2\lambda_n)}\right)^2\right)^{\!1/2}\\ 
&\leq \sup_{n\in\mathbb{N}} \{\lambda_n^{-\nu}|B\varphi_n|\}\, 
\|f\|_{\dot{H}^{2\mu+2\nu}(\Omega)}
s^{-(2-\alpha_N+2/p)} \frac{c^2+\sum_{k=1}^N b_{k}(\lambda_1)^{\beta_{k}-1}}{
p^{1/p}(\tfrac{p}{p-1}c^2)^{\frac{p-1}{p}}}
\left(\sum_{n=1}^\infty\lambda_n^{-2(\mu-1/p)}\right)^{\!1/2} 
\end{aligned}
\]
where we have used Young's inequality to estimate
\[
s^2+c^2\lambda_n \geq 
s^{2/p}p^{1/p}(\tfrac{p}{p-1}c^2\lambda_n)^{\frac{p-1}{p}}.
\]
Since $\lambda_n \sim n^{2/d}$, the sum converges for $2(\mu-\frac{1}{p})>\frac{d}{2}$. 
On the other hand, to get an $O(s^{-(2+\gamma_J+\epsilon)}) $ estimate for $D$ with $\epsilon>0$, we need $2-\alpha_N+2/p>2+\gamma_J$ and therefore choose
\begin{equation}\label{eqn:mu}
\mu>\frac{d}{4}+\frac{\max\{\alpha_N+\gamma_J,\tilde{\alpha}_{\tilde{N}}+\tilde{\gamma}_{\tilde{J}}\}}{2}\,, \quad p:= \frac{2}{\max\{\alpha_N+\gamma_J,\tilde{\alpha}_{\tilde{N}}+\tilde{\gamma}_{\tilde{J}}\}}
\end{equation} 
(taking into account the tilde version of the above estimate as well).
This yields
\[
\begin{aligned}
&\frac{\sum_{n=1}^\infty\langle f_{i}\varphi_n\rangle B_i\varphi_n}{
s^2+c^2\lambda_1+\sum_{j=1}^J d_j s^{\gamma_j+2}+\sum_{k=1}^{N} s^{\alpha_k} b_{k}(\lambda_1)^{\beta_{k}}}\\
&\qquad=
\frac{\sum_{n=1}^\infty\langle \tilde{f}_{i}\tilde{\varphi}_n\rangle \tilde{B}_i\tilde{\varphi}_n}{
s^2+\tilde{c}^2\tilde{\lambda}_1+\sum_{j=1}^{\tilde{J}} \tilde{d}_j s^{\tilde{\gamma}_j+2}+\sum_{k=1}^{\tilde{N}} s^{\tilde{\alpha}_k} \tilde{b}_{k}(\tilde{\lambda}_1)^{\tilde{\beta}_{k}}}
+O(s^{-(2+\max\{\gamma_J,\tilde{\gamma}_{\tilde{J}}\}+\epsilon)}), \quad s\in[1,\infty) 
\end{aligned}
\]
with $\epsilon = \mu-\frac{d}{4}-\frac{\max\{\alpha_N+\gamma_J,\tilde{\alpha}_{\tilde{N}}+\tilde{\gamma}_{\tilde{J}}\}}{2}>0$. 
Thus, our assumptions on $\mu$, $\nu$ are 
\begin{equation}\label{eqn:munu}
\mu>\frac{d}{4}+\frac{\max\{\alpha_N+\gamma_J,\tilde{\alpha}_{\tilde{N}}+\tilde{\gamma}_{\tilde{J}}\}}{2},
\quad
\sup_{n\in\mathbb{N}} \{\lambda_n^{-\nu}|B\varphi_n|\}, \ \sup_{n\in\mathbb{N}} \{\tilde{\lambda}_n^{-\nu}|B\tilde{\varphi}_n|\} <\infty.
\end{equation}

In case of rational powers 
\begin{equation}\label{eqn:rationalpowers}
\begin{aligned}
&\alpha_k=\frac{p_k}{q_k}, \ \tilde{\alpha}_k=\frac{\tilde{p}_k}{\tilde{q}_k}, \ \gamma_j=\frac{n_j}{m_j}, \  \tilde{\gamma}_j=\frac{\tilde{n}_j}{\tilde{m}_j}, \\
&\bar{q}=\mbox{lcm}(q_1,\ldots,q_{N},\tilde{q}_1,\ldots,\tilde{q}_{\tilde{N}},m_1,\ldots,m_{J},\tilde{m}_1,\ldots,\tilde{m}_{\tilde{J}}),
\end{aligned}
\end{equation} 
setting $z=s^{1/\bar{q}}$ we read \eqref{eqn:BaBb_singlemodes} as an equality between two rational functions of $z$, whose coefficients and powers therefore need to coincide (provided $\{s^{1/\bar{q}}\, : \, s\in\mathbb{C}_\theta\}\supseteq [\underline{z},\infty)$ for some $\underline{z}\geq0$).
In the general case of real powers one can use the induction proof from \cite[Theorem 1.1]{JinKian:2021} to obtain the same uniqueness.
This yields
\[
\begin{aligned}
&(I)\qquad&&N=\tilde{N},\qquad J=\tilde{J}\,,\\
&(II)&&c^2\lambda_1=\tilde{c}^2\tilde{\lambda}_1\,, \\
&(III)&&\alpha_k=\tilde{\alpha}_k\,,\ \gamma_j=\tilde{\gamma}_j\,, &&k=1,\ldots N\,, \ j=1,\ldots J\\
&(IV)&&b_k(\lambda_1)^{\beta_{k}}
=\tilde{b}_k(\tilde{\lambda}_1)^{\tilde{\beta}_{k}}\,, \ d_j=\tilde{d}_j\,, \quad 
&& \ k=1,\ldots N\,, \ j=1,\ldots J
\end{aligned}
\]
From this it becomes clear that we get uniqueness of 
\[
N, \ J, \ c, \ b_{k}, \ \alpha_k, \ d_j, \ \gamma_j,
\]
provided all $\beta_k=\tilde{\beta}_k$ are known and $\lambda_1=\tilde{\lambda}_1$ (but still possibly unknown) since we can then divide by $\lambda_1$ and $\lambda_1^{\beta_k}$ in (II) and (IV), respectively.
However, there seems to be no chance to simultaneously obtain $b_k$ and $\beta_k$ from (IV). Therefore, we will consider a setting with two different excitations in Section~\ref{sec:singlemode} for this purpose.

The same derivation can be made for the variable $c=c(x)$ model \eqref{eqn:general_c} in place of \eqref{eqn:general}, with 
$c$, $\mathcal{A}$, $L^2(\Omega)$, $\|\varphi\|_{L^2}=1$ replaced by 
$1$, $\mathcal{A}_c$, $L^2_{1/c^2}(\Omega)$, $\|\varphi\|_{L^2_{1/c^2}}=1$. 
From (II) with $c$, $\tilde{c}$ replaced by unity, we get $\lambda_1=\tilde{\lambda}_1$ so that we do not need to assume this. 

\begin{theorem}\label{thm:uniqueness_smoothdata}
Let $f,\tilde{f}\in \dot{H}^{2\mu+2\nu}(\Omega)$ for $\mu,\nu$ sufficiently large so that \eqref{eqn:munu} holds.
 
Then
\[
B u(t) = \tilde B \tilde u(t)\, \quad t\in(0,T)
\]
\begin{enumerate}
\item[(i)] 
for the solutions $u$, $\tilde{u}$ of \eqref{eqn:general}
with vanishing initial data and known equal $\beta_k=\tilde{\beta}_k$, possibly unknown but equal $\lambda_1=\tilde{\lambda}_1$, possibly unknown (rest of) $\mathcal{A}$, $\tilde{\mathcal{A}}$, and possibly unknown $f,\tilde{f}$
implies
\[
N=\tilde{N}, \ J=\tilde{J}, \ c=\tilde{c}, \ b_k=\tilde{b}_k, \ \alpha_k=\tilde{\alpha}_k, \ d_j=\tilde{d}_j, \ \gamma_j=\tilde{\gamma}_j,\quad k\in\{1,\ldots,N\}, j\in\{1,\ldots,J\}\,. 
\]
\item[(ii)]
for the solutions $u$, $\tilde{u}$ of \eqref{eqn:general_c}
with vanishing initial data and known equal $\beta_k=\tilde{\beta}_k$, and possibly unknown $\mathcal{A}_c$, $\tilde{\mathcal{A}}_{\tilde{c}}$
implies
\[
N=\tilde{N}, \ J=\tilde{J}, \ b_k=\tilde{b}_k, \ \alpha_k=\tilde{\alpha}_k, \ d_j=\tilde{d}_j, \ \gamma_j=\tilde{\gamma}_j,\quad k\in\{1,\ldots,N\}, j\in\{1,\ldots,J\}\,. 
\]
\end{enumerate}
\end{theorem}

\begin{remark}\label{rem:uniqueness_u0smooth}
This uniqueness approach also extends to different excitation combinations such as 
$u_0\not=0$, $u_j=0$, $j=1,2$, $f=0$; 
$u_1\not=0$, $u_j=0$, $j=0,2$, $f=0$; 
$u_2\not=0$, $u_j=0$, $j=0,1$, $f=0$; (the latter only in case $\gamma_J>0$) 

Consider for example the case $u_0\not=0$, $u_j=0$, $j=1,2$, $f=0$ that will be relevant for additionally recovering $u_0$, cf. Remark~\ref{rem:uniqueness_smoothdata_and_u1_u0} and which is the one that differs most from the setting above in view of the fact that the numerator of the resolvent solution depends on $s$ and $\lambda_n$.
The crucial estimate on the error that is made by replacing $\lambda_n$ by $\lambda_1$ then becomes, for $s\in\mathbb{R}$, $s\geq1$, 
\revision{and using the decomposition $\omega(s,\lambda)=s^2+\omega_d(s)+\omega_b(s,\lambda)$ with $\omega_d(s)=\sum_{j=1}^J d_j s^{\gamma_j+2}$, $\omega_b(s,\lambda)=\sum_{k=1}^N s^{\alpha_k} b_{k}\lambda^{\beta_{k}}$
\[
\begin{aligned}
D:=&\left|
\sum_{n=1}^\infty
\frac{\langle u_0,\varphi_n\rangle B\varphi_n(\omega(s,\lambda_1)-c^2\lambda_1)}{s\omega(s,\lambda_1)}-
\sum_{n=1}^\infty\frac{\langle u_0,\varphi_n\rangle B\varphi_n(\omega(s,\lambda_n)-c^2\lambda_n)}{s\omega(s,\lambda_n)}
\right|\\
&=\left|\sum_{n=1}^\infty
\frac{\langle u_0,\varphi_n\rangle B\varphi_n\,c^2\Bigl(\lambda_n(s^2+\omega_d(s)+\omega_b(s,\lambda_1))-\lambda_1(s^2+\omega_d(s)+\omega_b(s,\lambda_n))\Bigr)}{s\,\omega(s,\lambda_n)\, \omega(s,\lambda_1)}\right|\\
\\
&\leq\sum_{n=1}^\infty
\frac{|\langle u_0,\varphi_n\rangle B\varphi_n|\,c^2\Bigl((\lambda_n-\lambda_1)(s^2+\omega_d(s))+\lambda_n\omega_b(s,\lambda_1)-\lambda_1\omega_b(s,\lambda_n)\Bigr)}{s\,\omega(s,\lambda_n)\, \omega(s,\lambda_1)}\\
&\leq \sup_{n\in\mathbb{N}} \{\lambda_n^{-\nu}|B\varphi_n|\}\, 
\|u_0\|_{\dot{H}^{2\mu+2\nu}(\Omega)}
\left(\sum_{n=1}^\infty\left(
\frac{(s^2+\bar{d} s^{2+\gamma_J} +\bar{b}s^{\alpha_N})c^2\lambda_n}{
\lambda_n^\mu s^3(s^2+c^2\lambda_n)}\right)^2\right)^{1/2}
\end{aligned}
\]
for $\bar{d}=\sum_{j=1}^J d_j$, $\bar{b}=\sum_{k=1}^N b_{k}\max\{1,\lambda_1\}$,
where we have used $(\lambda_n-\lambda_1)(s^2+\omega_d(s))\geq0$ and $\lambda_n\omega_b(s,\lambda_1)-\lambda_1\omega_b(s,\lambda_n)
=\sum_{k=1}^N s^{\alpha_k} b_{k}\bigl(\lambda_n^{\beta_{k}}\lambda_1^{\beta_{k}}(\lambda_n^{1-\beta_{k}}-\lambda_1^{1-\beta_{k}})\bigr)\geq0$.
}
\\
To achieve an $O(s^{-(2+\max\{\gamma_J,\tilde{\gamma}_{\tilde{J}}\}+\epsilon)})$ estimate, instead of \eqref{eqn:munu} we thus assume
\begin{equation}\label{eqn:munu_u0}
\mu>\frac{d}{4}+1,\quad \max\{\gamma_J,\tilde{\gamma}_{\tilde{J}}\}<\frac12\, \quad
\quad
\sup_{n\in\mathbb{N}} \{\lambda_n^{-\nu}|B\varphi_n|\}, \ \sup_{n\in\mathbb{N}} \{\tilde{\lambda}_n^{-\nu}|B\tilde{\varphi}_n|\} <\infty.
\end{equation}
to recover Theorem~\ref{thm:uniqueness_smoothdata} with $f$, $\tilde{f}$ replaced by $u_0$, $\tilde{u}_0$.
\end{remark}

\begin{remark}
Note that in view of the fact that $\mathcal{A}$, $\mathcal{A}_c$ do not need to be known (up to the fact that a single eigenvalue is supposed to coincide), we are able to identify all relevant information on the damping model in an unknown medium.
\end{remark}

\subsection{Single mode excitation}\label{sec:singlemode}
Assume that we know at least some of the eigenfunctions and can use them as excitations $f_{i}=\varphi_{n_i}$, $\tilde{f}_{i}=\tilde{\varphi}_{\tilde{n}_i}$, $i=1,\ldots, I$, so that \eqref{eqn:BaBb} becomes
\begin{equation}\label{eqn:BaBb_singlemodes}
\begin{aligned}
&\frac{\langle f_{i},\varphi_{n_i}\rangle B_i\varphi_{n_i}}{
s^2+c^2\lambda_{n_i}+\sum_{j=1}^J d_j s^{\gamma_j+2}+\sum_{k=1}^N s^{\alpha_k} b_{k}(\lambda_{n_i})^{\beta_{k}}}\\
&\qquad=
\frac{\langle {\tilde f}_{i},{\tilde \varphi}_{\tilde{n}_i}\rangle {\tilde B}_i {\tilde \varphi}_{\tilde{n}_i}}{
s^2+{\tilde c}^2 {\tilde \lambda}_{\tilde{n}_i}+\sum_{j=1}^{\tilde{J}} \tilde{d}_j s^{\tilde{\gamma}_j+2}+\sum_{k=1}^{\tilde N} s^{{\tilde \alpha}_k} {\tilde b}_{k}({\tilde \lambda}_{\tilde{n}_i})^{{\tilde \beta}_{k}}}, \quad i=1,\ldots I\,, \quad s\in\mathbb{C}_\theta\,. \end{aligned}
\end{equation}

We can use the rational function or induction argument as above to compare powers of $s$ and conclude the following:
\[
\begin{aligned}
&(I)\qquad&&N=\tilde{N},\qquad J=\tilde{J}\\
&(II)&&c^2\lambda_{n_i}=\tilde{c}^2\tilde{\lambda}_{\tilde{n}_i}\,, \quad&& i=1,\ldots,I\\
&(III)&&\alpha_k=\tilde{\alpha}_k\,,\ \gamma_j=\tilde{\gamma}_j\,, &&k=1,\ldots N\,, \ j=1,\ldots J\\
&(IV)&&b_k(\lambda_{n_i})^{\beta_{k}}
=\tilde{b}_k(\tilde{\lambda}_{\tilde{n}_i})^{\tilde{\beta}_{k}}\,, \ d_j=\tilde{d}_j\,, \quad 
&&i=1,\ldots,I\,, \ k=1,\ldots N\,, \ j=1,\ldots J
\end{aligned}
\]
To extract $b_k$ and $\beta_k$, which can be done separately for each $k$, (skipping the subscript $k$) we use $I=2$ excitations, assume that the corresponding eigenvalues of $\mathcal{A}$ and $\tilde{\mathcal{A}}$ are known and equal and define 
\[
F(b,\beta)=
(b(\lambda_{n_1})^\beta,b(\lambda_{n_2})^\beta)^T\,.
\]
One easily sees that the $2\times2$ matrix $F'(b,\beta)$ is regular for $\lambda_{n_1}\neq\lambda_{n_2}$ and $b\not=0$ and thus by the Inverse Function Theorem $F$ is injective. Thus from (I)-(IV) we obtain all the constants in \eqref{eqn:model}.

\begin{theorem}\label{thm:uniqueness_singlemode}
Let $f_1=\varphi_{n_1}$, $f_2=\varphi_{n_2}$, $\tilde{f}_1=\tilde{\varphi}_{\tilde{n}_1}$, $\tilde{f}_2=\tilde{\varphi}_{\tilde{n}_2}$ for some $n_1\not=n_2$, $\tilde{n}_1\not=\tilde{n}_2$ $\in\mathbb{N}$.
Then 
\[
B_i u_i(t) = \tilde B_i \tilde u_i(t), \quad t\in(0,T),\quad i=1,2
\]
for the solutions $u_i$, $\tilde{u}_i$ of \eqref{eqn:general}
with vanishing initial data, known and equal $\lambda_i=\tilde{\lambda}_i$, $i=1,2$ and possibly unknown (rest of) $\mathcal{A}$, $\tilde{\mathcal{A}}$
implies
\[
N=\tilde{N}, \ J=\tilde{J}, \ c=\tilde{c}, \ b_k=\tilde{b}_k, \ \alpha_k=\tilde{\alpha}_k, \ \beta_{k}=\tilde{\beta}_k, \ d_j=\tilde{d}_j, \ \gamma_j=\tilde{\gamma}_j,\quad k\in\{1,\ldots,N\}, j\in\{1,\ldots,J\}\,. 
\]
\end{theorem}

\begin{remark}
Scaling with unity in the excitation has been chosen here just for simplicity of exposition. Clearly also multiples of eigenfunctions can be used.

Again, uniqueness can be shown analogously with eigenfunction excitation by initial data rather than by sources.
\end{remark}

We briefly consider the additional recovery of a spatially varying wave speed $c$ and for this purpose look at  the {\sc ch} model \eqref{eqn:CH_c} as an example.
With $\mathcal{A}_c=-c(x)^2\triangle$, (I)--(IV) (with $N=1$, $c=1$, $\beta=1$, $b_{k}=b$) yields
\[
\alpha=\tilde{\alpha}\,, \quad 
\lambda_{n_i}=\tilde{\lambda}_{\tilde{n}_i}\,, \quad 
b\,\lambda_{n_i}=\tilde{b}\,\tilde{\lambda}_{\tilde{n}_i}\,, \quad i=1,\ldots,I\,.
\]
Using inverse Sturm-Liouville theory, one can recover the spatially varying coefficient $c(x)$ in one space dimension from knowledge of all eigenvalues (actually just two sets of eigenvalues corresponding to different boundary conditions),
see, \cite{RundellSacks:1992b,CCPR:1997}.
However, in order to do so via this single mode exciation approach, we would need infinitely many measurements
$\{n_i\, : i=1,\ldots I\}=\mathbb{N}$.

The goal of obtaining knowledge on the eigenvalues can be much better achieved by going back to \eqref{eqn:BaBb} and considering equality of poles (and residues), which according to 
\cite[Section 4]{KaltenbacherRundell:2021b}
gives us $\lambda_n$ and $b\lambda_n^\beta$ (note that we have set $c=1$ since $c(x)$ is contained in $\mathcal{A}_c$) from a single excitation $u_1$ with nonvanishing Fourier coefficients. We will now consider this approach for the more general setting \eqref{eqn:general_2ndorder} in Section~\ref{sec:poles}.

\subsection{Poles and residues}\label{sec:poles}
In this section we focus on the second order in time case
\eqref{eqn:general_2ndorder} and achieve separation of summands in the
eigenfunction expansion by using poles (and residues) of $\hat{h}(s)$.
For this purpose we have to prove that each summand $\lambda_n$ gives rise to at
least one pole (in case of {\sc ch} we know there are two complex conjugates) and
that these poles differ for different $\lambda_n$.
This is also known for {\sc ch} and {\sc fz}, see 
\cite[Section 4]{KaltenbacherRundell:2021b}.
Indeed the results from 
\cite[Lemma 4.1, Remark 4.1]{KaltenbacherRundell:2021b}
on {\sc ch} carry over to the model \eqref{eqn:general_2ndorder} in the setting
\begin{equation}\label{eqn:samebeta}
\beta_1=\cdots=\beta_k=\beta, \quad b_1,\cdots,b_k\geq0 
\end{equation}
that is, \eqref{eqn:general_2ndorder} with a single $\beta$ and dissipative behaviour (as is natural to demand from a damping model). 
Existence of at least one pole for each $\lambda_n$ can even be shown for the most
general model \eqref{eqn:mostgeneral_2ndorder}. 
\begin{lemma}\label{lem:what}
For each $\lambda_n$ there exists at least one root of 
$\omega^{\text{mg}}_n(s)=s^2+c^2\lambda_n+\sum_{k=1}^N s^{\alpha_k}\sum_{\ell=1}^{M_k} c_{k\ell}\lambda_n^{\beta_{k\ell}}$, that is, a pole of the Laplace transformed relaxation solution $\hat{w}^{\text{mg}}_{fn}$, $\hat{w}^{\text{mg}}_{0n}$, or $\hat{w}^{\text{mg}}_{1n}$.
Moreover, the poles are single.
\\
In case of \eqref{eqn:samebeta}, the poles of $\hat{w}^{\text{g}}_{fn}$, $\hat{w}^{\text{g}}_{0n}$, $\hat{w}^{\text{g}}_{1n}$, which are the roots of $\omega^{\text{g}}_n(s)=s^2+c^2\lambda_n+\sum_{k=1}^N b_k s^{\alpha_k} \lambda_n^\beta$, lie in the left half plane and differ for different $\lambda_n$.  
\end{lemma}
\begin{proof}
Let $f_n(z) = z^2 + c^2\lambda_n$, $g_n(z) = \sum_{k=1}^N z^{\alpha_k}\sum_{\ell=1}^{M_k} c_{k\ell}\lambda_n^{\beta_{k\ell}}$.
Let $C_R$ be the circle radius $R$, centre at the origin.
Then, due to $\alpha_k\leq1$, for a sufficiently large $R>c^2\lambda_n$, the estimate $|g_n(z)| < |f_n(z)|$ holds on $C_R$ and so Rouch\'{e}'s theorem shows that
$f_n(z)$ and $\omega^{\text{mg}}_n(z)=(f_n+g_n)(z)$ have the same number of roots, counted with multiplicity, within $C_R$.
For $f_n$ these are only at $z=\pm i\sqrt{\lambda_n}c$
and so $\omega^{\text{mg}}_n$ has 
precisely one single root in the third and in the fourth quadrant, respectively.
\\
The fact that the poles of $\hat{w}^{\text{g}}_{fn}$, $\hat{w}^{\text{g}}_{0n}$, $\hat{w}^{\text{g}}_{1n}$ lie in the left half plane in case of \eqref{eqn:samebeta} follows by an energy argument similar to the one in the proof of Lemma 11.3 fde-lect.
Suppose now that $\omega^{\text{g}}$ has a root at $p=r e^{i\theta}$,
where $\pi/2\leq\theta<\pi$, for both $\lambda_n$ and $\lambda_m$. Then for $p=r e^{i\theta}$, subtracting the equations $\omega^{\text{g}}_n(p)=0$ and $\omega^{\text{g}}_m(p)=0$ we obtain
\begin{equation}\label{eqn:lambdan-lambdam}
\frac{c^2(\lambda_n - \lambda_m)}{\lambda_n^\beta-\lambda_m^\beta} = - \sum_{k=1}^N b_k p^{\alpha_k} = \frac{p^2+c^2\lambda_n}{\lambda_n^\beta}
\end{equation}
thus
$$
\frac{p^2}{c^2}
=\frac{\lambda_n^\beta(\lambda_n - \lambda_m)}{\lambda_n^\beta-\lambda_m^\beta}-\lambda_n
=\frac{\lambda_n^\beta\lambda_m^\beta(\lambda_n^{1-\beta}-\lambda_m^{1-\beta})}{\lambda_n^\beta-\lambda_m^\beta}
$$
Now if $\lambda_n \not=\lambda_m$ then the right hand side is positive and real
and so $2\theta =\pi$.
This means that $\sum_{k=1}^N b_k s^{\alpha_k}= \sum_{k=1}^N b_k r^{\alpha_k\pi/2}$ (where all $b_k$ have the same sign) has nonvanishing imaginary part, a contradiction to the fact that the left hand side of \eqref{eqn:lambdan-lambdam} is real.
\end{proof}
Since we will separate the summands by taking residues on both sides of the identity \eqref{eqn:BaBb}, the following identity is also crucial.
\begin{lemma}\label{lem:residues}
For each $\lambda_n$, $n\in\mathbb{N}$, the residues of $\hat{w}^{\text{g}}_{fn}$, $\hat{w}^{\text{g}}_{0n}$, $\hat{w}^{\text{g}}_{1n}$ do not vanish.
\end{lemma}
\begin{proof}
By l'Hospital's rule and due to the fact that the poles $p_n$ are single, we get for $\omega(s)=\omega^{\text{mg}}_{fn}(s)=\omega^{\text{mg}}_{1n}(s)$
\[
\begin{aligned}
\mbox{Res}(\hat{w}^{\text{mg}}_n;p_n)&=\lim_{s\to p_n}\frac{(s-p_n)}{\omega(s)}
=\lim_{s\to p_n}\frac{1}{\omega'(s)}\not=0\,.
\end{aligned}
\]
Moreover, with $\omega(s)=\omega^{\text{mg}}_{0n}(s)$, $d(s)=\frac{\omega(s)-c\lambda^2}{s}$,
\[
\begin{aligned}
\mbox{Res}(\hat{w}^{\text{mg}}_n;p_n)&=\lim_{s\to p_n}\frac{(s-p_n)d(s)}{\omega(s)}
=\lim_{s\to p_n}\frac{d(s)+(s-p_n)d'(s)}{\omega'(s)}=\frac{d(p_n)}{\omega'(p_n)}\\
&=\frac{\omega(p_n)-c\lambda_n^2}{p_n\omega'(p_n)}=-\frac{c\lambda_n^2}{p_n\omega'(p_n)}\not=0\,.
\end{aligned}
\]
\end{proof}

As a consequence of Lemma~\ref{lem:what}, \ref{lem:residues}, by taking residues at the poles $p_n$ of the Laplace transformed time trace data $\hat{h}(s)$, 
we obtain from \eqref{eqn:BaBb} 
\begin{equation}\label{eqn:BaBb_residues}
\begin{aligned}
\mbox{Res}(\hat{h};p_n)=&\frac{1}{
2p_n+\sum_{k=1}^N \alpha_k b_k p_n^{\alpha_k-1}\, (\lambda_{n})^{\beta_k}}
\sum_{k\in K^{\lambda_n}}\langle f,\varphi_{n,k}\rangle B\varphi_{n,k}
\\&\qquad
=\frac{1}{
2p_n+\sum_{k=1}^{\tilde N} {\tilde \alpha}_k \tilde{b}_k p_n^{{\tilde \alpha}_k-1}({\tilde \lambda}_{n})^{\tilde{\beta}_k}}
\sum_{k\in K^{\tilde{\lambda}_n}} \langle {\tilde f},{\tilde \varphi}_{n,k}\rangle {\tilde B} {\tilde \varphi}_{n,k}
, \\ 
&i=1,\ldots I\,, \quad n\in\mathbb{N}\,. 
\end{aligned}
\end{equation}
\revision{where the cardinality of the index set $K^{\lambda_n}$ is equal to the multiplicity of the eigenvalue $\lambda_n$.}

In order to achieve nonzero numerators in \eqref{eqn:BaBb}, \eqref{eqn:BaBb_singlemodes} and \eqref{eqn:BaBb_residues}, we assume
\begin{equation}\label{eqn:Bnonzero}
B\varphi_{n}\not=0, \quad \tilde{B}\tilde{\varphi}_{n} \not=0\quad\mbox{ for all }n\in\mathbb{N},
\end{equation}
and
\begin{equation}\label{eqn:fmodesnonzero}
\langle f,\varphi_n\rangle\not=0, \quad \langle \tilde{f},\tilde{\varphi}_n\rangle\not=0\quad \mbox{ for all }n\in\mathbb{N},
\end{equation} 
\revision{see remark~\ref{rem:nonzerocoeff} below.}
\Margin{Ref 2 (iii)}

Assume that, 
\begin{equation}\label{eqn:fBequal}
\sum_{k\in K^{\lambda_n}}\langle f,\varphi_{n,k}\rangle B\varphi_{n,k}
=\sum_{k\in K^{\tilde{\lambda}_n}} \langle {\tilde f},{\tilde \varphi}_{n,k}\rangle {\tilde B} {\tilde \varphi}_{n,k}\not=0\,, \quad 
\mbox{ for all }n\in\mathbb{N}.
\end{equation}
(e.g. by assuming $f=\tilde{f}$, $B=\tilde{B}$, $\mathcal{A}=\tilde{\mathcal{A}}$). 
Then from the reciprocals of \eqref{eqn:BaBb_residues} we get, after division by 
$\sum_{k\in K^{\lambda_n}}\langle f,\varphi_{n,k}\rangle B\varphi_{n,k}$ and subtracting $2p_n$ on both sides
\begin{equation}\label{eqn:BaBb_residues_recip}
\sum_{k=1}^N \alpha_k b_k p_n^{\alpha_k-1}\, (\lambda_{n})^{\beta_k}
=\sum_{k=1}^{\tilde N} {\tilde \alpha}_k \tilde{b}_k p_n^{{\tilde \alpha}_k-1}(\tilde{\lambda_{n}})^{\tilde{\beta}_k}.
\end{equation}
In case $\beta=\beta_k=\tilde{\beta}_k=\tilde{\beta}$, $\lambda_n=\tilde{\lambda}_n$, this can be divided by $(\lambda_{n})^{\beta}$ and simplifies to 
\[
\sum_{k=1}^N \alpha_k b_k p_n^{\alpha_k-1}
=\sum_{k=1}^{\tilde N} {\tilde \alpha}_k \tilde{b}_k p_n^{{\tilde \alpha}_k-1}.
\]
that is, in case of rational powers \eqref{eqn:rationalpowers}
an equality between two polynomials of $z=s^{1/\bar{q}}$ at infinitely many different (according to the second part of Lemma~\ref{lem:what}) points $p_n^{1/\bar{q}}$. 
This yields
\[
N=\tilde{N}\,,\quad \alpha_k=\tilde{\alpha}_k\,, \quad b_k=\tilde{b}_k\,,\quad  k=1,\ldots N.
\]
Thus we recover the results of Theorem~\ref{thm:uniqueness_smoothdata} under essentially more restrictive assumptions 
-- however we can drop the smoothness assumption on the excitation $f$, so that $f=\delta$ is admissible.
Thus the use of poles and residues has no clear advantage over the smooth excitation approach as long as we only aim at recovering the damping term constants. However, it allows to gain additional information on spatially vaying coefficients $c=c(x)$ or $f=f(x)$ or $u_0=u_0(x)$.

\subsubsection{Additional recovery of $c(x)$}

First of all, just using knowledge of the pole locations of the data,
Lemma~\ref{lem:what} allows to obtain the following corollary of
Theorem~\ref{thm:uniqueness_smoothdata} (ii) on uniqueness of the numbers $N$, $J$,
coefficients $b_k$, orders $\alpha_k$, as well as the function $c(x)$ in one space
dimension. Note that the assumption $f=\delta$ is not compatible with the
smoothness assumptions from Theorem \ref{thm:uniqueness_smoothdata}, thus we cannot
apply the strategy from \cite{KaltenbacherRundell:2021b} of obtaining, besides the
eigenvalues, also the values $\varphi(x_0)$ and therewith the norming constants of
the eigenfunctions.
We thus apply a different result from Sturm-Liouville theory,
according to which a pair of eigenfunction sequences corresponding to two different
boundary impedances  suffices to recover the potential $q$ of the differential
operator $-\triangle + q$, see \cite{RundellSacks:1992a,CCPR:1997}.
Using the transformation of variables 
\[
\xi(x)=\int_0^xc(s)^{-1}\, ds, \quad \eta(\xi)=c^{-1/2}(x) y(x), \quad q(\xi)= \frac14(c'(x)^2-2c(x)c''(x)),
\]
so that
\[
-c^2(x)y_{xx}(x)=\lambda y(x), \ x\in (0,L) \ \Leftrightarrow \
-\eta_{\xi\xi}(\xi)+q(\xi)\eta(\xi)=\lambda\eta(\xi), \ \xi\in (0,\int_0^Lc(s)^{-1}\, ds) 
\]
this transfers to the recovery of the sound speed $c(x)$.

Since these results hold in one space dimension only, in this subsection we consider 
\begin{equation}\label{eqn:general_c_oned}
\begin{aligned}
&u_{tt}-c(x)\triangle u
\sum_{k=1}^N b_k \partial_t^{\alpha_k} (-c(x)\triangle u)^{\beta_k}u=\sigma(t)f\quad && x\in (0,L), \ t\in(0,T)\\
&u'(0,t)-H_0u(0,t)=0\,,\quad u'(L,t)+H_Lu(L,t)=0\quad && t\in(0,T)\\
&u(x,0)=u_0(x)\,, \quad u_t(x,0)=u_1(x)\, \quad u_{tt}(x,0)=u_2(x)\quad &&x\in(0,L)
\end{aligned}
\end{equation}
for some $L>0$ and some nonnegative impedance constants $H_0$, $H_L$.
For two different right hand boundary impedance values $H_{L,i}$, $i=1,2$ we denote by $u^i$ the solution of \eqref{eqn:general_c_oned} with $H_L=H_{L,i}$ and by 
$(\varphi_{n,i}$, $\lambda_{n,i})_{n\in\mathbb{N}}$ the eigenpairs of $-c(x)\triangle$ with boundary conditions $\varphi_{n,i}'(0)-H_0\varphi_{n,i}(0)=0$, $\varphi_{n,i}'(L)+H_{L,i}\varphi_{n,i}(L)=0$, 
and weigthed $L^2$ normalization $\int_0^L \frac{1}{c(x)^2}\varphi_{n,i}(x)^2\, dx = 1$, $i=1,2$.

For constructing sufficiently smooth exciations with nonvanishing coefficients with respect to the eigenfunction basis, 
we point to Remark~\ref{rem:nonzerocoeff} below.

Moroever, we assume that the observations satisfy 
\begin{equation}\label{eqn:Binonzero}
B_i\varphi_{n,i}\not=0, \quad \tilde{B}_i\tilde{\varphi}_{n,i} \not=0\mbox{ for all }n\in\mathbb{N}, \quad i=1,2, 
\end{equation}
\revision{which is, e.g., the case with the choice   
$B_i v:= v(x_0)$, $x_0\in \{0,L\}$,
see remark~\ref{rem:nonzerocoeff} below.}
\Margin{Ref 2 (iii)}

Analogous notation is used for \eqref{eqn:general_c_oned} with $c$, 
$b_k$, $\alpha_k$, $\beta_k$, $f$, $H_0$, $H_L$ replaced by their corresponding tilde counterparts.

From Theorem~\ref{thm:uniqueness_smoothdata} and uniquenss of the eigenvalues from the poles we obtain the following result on combined identification of the coefficients in the damping model and the space-dependent sound speed.
\begin{corollary}\label{cor:uniqueness_smoothdata_and_c}
Let $f,\tilde{f}\in \dot{H}^{2\mu+2\nu}(0,L)$ for $\mu,\nu$ satisfying \eqref{eqn:munu} and let
\eqref{eqn:fmodesnonzero}, \eqref{eqn:Binonzero} hold.

Then
\[
B_i u_i(t) = \tilde B_i \tilde u_i(t)\, \quad t\in(0,T), \quad i=1,2
\]
for the solutions $u_i$, $\tilde{u}_i$ of \eqref{eqn:general_c_oned}
with vanishing initial data, $H_L=H_{L,i}$, $\tilde{H}_L=\tilde{H}_{L,i}$, known equal $\beta_k=\tilde{\beta}_k$, and unknown $c$, $\tilde{c}$ $\in C^2(0,L)$,
implies
\[
\begin{aligned}
&N=\tilde{N}, \ 
b_k=\tilde{b}_k, \ \alpha_k=\tilde{\alpha}_k, \ 
\quad k\in\{1,\ldots,N\}, 
\\ 
&c(x)=\tilde{c}(x), \quad x\in(0,L)\,. 
\end{aligned}
\]
\end{corollary}
The same result can be obtained without smoothness assumption on $f$ 
provided \eqref{eqn:fBequal}, and $\beta=\beta_k=\tilde{\beta}_k=\tilde{\beta}$ by a poles/residual argument, cf. \eqref{eqn:BaBb_residues_recip}.

\revision{
\begin{remark}\label{rem:nonzerocoeff}
Some comments on the verification of conditions \eqref{eqn:fmodesnonzero}, \eqref{eqn:Binonzero}, \eqref{eqn:fBequal} are in order, in particular in view of the fact that $c(x)$, $\tilde{c}(x)$ and therefore also the eigenfunctions are unknown.
\\
First of all, note that choosing $f,\tilde{f}=\delta_{x_0}$ formally satisfies this condition provided we can make sure that none of the eigenfunctions vanish at the point $x_0$. This is, e.g., the case if $x_0\in\{0,L\}$ is a boundary point for then in case $x_0=L$ the impedance conditions 
would imply 
$\varphi_n'(0)-H_0\varphi_n(0)=0$, $\varphi_n'(L)=0$, $\varphi_n(L)=0$, which by uniqueness of solutions to the second order ODE $-c(x)\varphi_n''(x)-\lambda_n\varphi_n(x)=0$ with these boundary conditions would imply $\varphi_n\equiv0$, a contradiction to $\varphi_n$ being an eigenfunction.
Likewise for $x_0=0$.
\\
Note that setting $B_i v:= v(x_0)$, $x_0\in \{0,L\}$, the same argument yields \eqref{eqn:Binonzero} as well as \eqref{eqn:fBequal}.
\\
However, clearly $f=\delta_{x_0}$ does not satisfy the regularity requirements from the first part of the corollary, thus we also indicate a possibility of constructing a smooth $f$ that satisfies \eqref{eqn:fmodesnonzero}.
For any $s>0$, $\epsilon>0$, the function 
$f(x)=c(x)^2\sum_{n=1}^\infty \lambda_n^{-s} n^{-(\frac12+\epsilon)} \varphi_n(x)$ lies in $\dot{H}^{s}$ and satisfies 
$\langle f,\varphi_k\rangle_{L^2_{1/c^2}}  \geq \lambda_n^{-s} n^{-(\frac12+\epsilon)} $ for all $k\in\mathbb{N}$.
Unfortunately this is only a result on existence of $f$, since the reconstruction requires knowledge of $c$.
\end{remark}
}
\Margin{Ref 2 (iii)}

\subsubsection{Additional recovery of $f(x)$ (or $u_0(x)$ or $u_1(x)$)}

In order to reconstruct the possibly unknown excitation $f$ (analogously for $u_1$, $u_0$, the latter relevant in PAT) we use measurements not only at finitely many points but on a surface $\Sigma$  
and make the linear independence assumption
\begin{equation}\label{eqn:ass_inj_Sigma_rem}
\left(\sum_{k\in K^{\lambda_n}} m_k \varphi_k(x) = 0 \ \mbox{ for all }x\in\Sigma\right)
\ \Longrightarrow \ \left(m_k=0 \mbox{ for all }k\in K^{\lambda_n}\right)
\quad n\in\mathbb{N}\,.
\end{equation}
Moreover, for the purpose of recovering $f$ (or $u_1$, or $u_0$) we assume $c(x)$ to be known and therefore $\varphi_{n}=\tilde{\varphi}_{n}$, $\lambda_{n}=\tilde{\lambda}_{n}$.
Thus in place of \eqref{eqn:BaBb_residues} we have
\begin{equation}\label{eqn:BaBb_residues_Sigma}
\begin{aligned}
&\frac{1}{
2p_n+\sum_{k=1}^N \alpha_k b_k p_n^{\alpha_k-1}\, (\lambda_{n})^\beta}
\sum_{k\in K^{\lambda_n}}\langle f,\varphi_{n,k}\rangle \varphi_{n,k}(x)
\\&\qquad
=\frac{1}{
2p_n+\sum_{k=1}^{\tilde N} {\tilde \alpha}_k \tilde{b}_k p_n^{{\tilde \alpha}_k-1}(\lambda_{n})^{\tilde{\beta}}}
\sum_{k\in K^{\tilde{\lambda}_n}}\langle \tilde{f},\varphi_{n,k}\rangle \varphi_{n,k}(x), 
\quad x\in\Sigma\,, \quad n\in\mathbb{N}\,. 
\end{aligned}
\end{equation}

Again, we can obtain a combined identification result from Theorem~\ref{thm:uniqueness_smoothdata}, uniquenss of the eigenvalues from the poles, and \eqref{eqn:BaBb_residues_Sigma}.

\begin{corollary}\label{cor:uniqueness_smoothdata_and_f}
Let $f,\tilde{f}\in \dot{H}^{2\mu+2\nu}(\Omega)$ for $\mu,\nu$ saftisfying \eqref{eqn:munu} and let  
\eqref{eqn:BaBb_residues_Sigma} hold.

Then
\[
u(x,t) = \tilde{u}(x,t)\, \quad x\in\Sigma, \ t\in(0,T), 
\]
for the solutions $u$, $\tilde{u}$ of \eqref{eqn:general_c}
with vanishing initial data known equal $\beta_k=\tilde{\beta}_k$, and unknown $f$, $\tilde{f}$ $\in L^2(\Omega)$,
implies
\[
\begin{aligned}
&N=\tilde{N}, \ 
b_k=\tilde{b}_k, \ \alpha_k=\tilde{\alpha}_k, \ 
\quad k\in\{1,\ldots,N\}, 
\\ 
&f(x)=\tilde{f}(x), \quad x\in\Omega\,. 
\end{aligned}
\]
\end{corollary}

\begin{remark}
The same result can be obtained without smoothness assumption on $f$ 
provided \eqref{eqn:fBequal}, and $\beta=\beta_k=\tilde{\beta}_k=\tilde{\beta}$ by a poles/residual argument, cf. \eqref{eqn:BaBb_residues_recip}.
\end{remark}

\begin{remark}\label{rem:uniqueness_smoothdata_and_u1_u0}
The result carries over to the setting $f=0$, $u_0=0$, where $u_1=u_1(x)\not=0$ is supposed to be recovered almost without changes.

In the setting $f=0$, $u_1=0$, of recovering $u_0=u_0(x)\not=0$ (as relevant in PAT), we can rely on Remark~\ref{rem:uniqueness_u0smooth} to conclude that Corollay~\ref{cor:uniqueness_smoothdata_and_f} remains valid with $f$, $\tilde{f}$, \eqref{eqn:munu} replaced by $u_0$, $\tilde{u}_0$, \eqref{eqn:munu_u0}.
\end{remark}

%% file: recon_meth.tex
\section{Reconstruction approaches}\label{sec:recon_meth}
Throughout this section we will a single mode excitation and therefore skip the subscript $n_i$ in $\varphi_{n_i}$, $\lambda_{n_i}$ and $i$ in $B_i$ to simplify notation). Moreover, we focus on the second order model \eqref{eqn:general_2ndorder}.

\subsection{Full time observations or analyticity: Laplace domain reconstruction}

Starting from the equation 
\[
\hat{h}(s)=\frac{\langle f,\varphi\rangle B\varphi}{
s^2+c^2\lambda+\sum_{k=1}^N s^{\alpha_k} b_{k}\lambda^{\beta_k}}
\begin{cases}
1 
&\mbox{ in case $\;u_0=0$, $\;u_1=\varphi$, $\;f=0$}\\
\hat{\sigma}(s)
&\mbox{ in case $\;u_0=0$, $\;u_1=0$, $\;f=\varphi$}\\
s+\sum_{k=1}^N s^{\alpha_k-1} b_{k}\lambda^{\beta_k}
&\mbox{ in case $\;u_0=\varphi$, $\;u_1=0$, $\;f=0$}
\end{cases}
\]
(cf. \eqref{eqn:BaBb_singlemodes} for the slightly simpler setting of \eqref{eqn:general}) and noting that taking the reciprocal and rearranging makes the problem somewhat less nonlinear and the derivative easier to compute, we consider 
\[
F(b_1,\ldots,b_N,\beta_1,\ldots,\beta_N,\alpha_1,\ldots,\alpha_N;s)
= G(s)\,, \quad s\in\{s_1,\ldots,s_M\}
\]
with 
\[
F(b_1,\ldots,b_N,\beta_1,\ldots,\beta_N,\alpha_1,\ldots,\alpha_N;s)
=\sum_{k=1}^N b_k \lambda^{\beta_k} s^{\alpha_k}
\]
and 
\[
G(s)=\begin{cases}
\frac{B\varphi}{\hat{h}(s)}-s^2-c^2\lambda
&\mbox{ in case $\;u_0=0$, $\;u_1=\varphi$, $\;f=0$}\\
\frac{B\varphi\, \hat{\sigma}(s)}{\hat{h}(s)}-s^2-c^2\lambda
&\mbox{ in case $\;u_0=0$, $\;u_1=0$, $\;f=\varphi$}\\
\frac{c^2\lambda}{1-\frac{B\varphi}{s\hat{h}(s)}}-s^2-c^2\lambda
&\mbox{ in case $\;u_0=\varphi$, $\;u_1=0$, $\;f=0$}
\end{cases}
\]
For applying Newton's method, the Jacobian is given by 
\[
\begin{aligned}
&\frac{\partial F}{\partial b_i}(b_1,\ldots,b_N,\beta_1,\ldots,\beta_N,\alpha_1,\ldots,\alpha_N;s):=
\lambda^{\beta_i} s^{\alpha_i}\\
&\frac{\partial F}{\partial \beta_i}(b_1,\ldots,b_N,\beta_1,\ldots,\beta_N,\alpha_1,\ldots,\alpha_N;s):=
b_i \ln \lambda\, \lambda^{\beta_i} s^{\alpha_i}\\
&\frac{\partial F}{\partial \alpha_i}(b_1,\ldots,b_N,\beta_1,\ldots,\beta_N,\alpha_1,\ldots,\alpha_N;s):=
b_i \lambda^{\beta_i} \ln s \, s^{\alpha_i}
\end{aligned}
\] 

Taking the logarithm makes the equation less nonlinear with respect to $\alpha$ and $\beta$ by considering (in case $N=1$)
\[
\begin{aligned}
F_1(b,\beta,\alpha)&:=
b \lambda^\beta s_1^\alpha
=\frac{B\varphi\, \hat{\sigma}(s_1)}{\hat{h}(s_1)}-
s_1^2-c^2\lambda\\
F_m(b,\beta,\alpha)&:=
\log b + \beta\log\lambda +\alpha \log s_m
=\log\left(\frac{B\varphi\, \hat{\sigma}(s_m)}{\hat{h}(s_m)}-
s_m^2-c^2\lambda\right)\,,\quad m=2,\ldots M.
\end{aligned}
\]
Then $F_1$ is linear with respect to $b$ and the $F_m$ are linear with respect to $\alpha$ and $\beta$. This helps to enlarge the convergence radius of Newton's method, as we experienced in our numerical tests.

\subsection{Large time asymptotics}

In case $u_0=\varphi$, $u_1=0$, $f=0$ we have 
\begin{equation}\label{eqn:hath}
\hat{h}(s) = \frac{s + \sum_{k=1}^N b_k\lambda^{\beta_k} s^{\alpha_k-1}}{s^2 + \sum_{k=1}^N b_k\lambda^{\beta_k} s^{\alpha_k} + c^2\lambda} B\varphi
\sim \frac{B\varphi}{c^2} \sum_{k=1}^N b_k\lambda^{\beta_k-1} s^{\alpha_k-1} \quad \mbox{ as }s\to0.
\end{equation}

\revision{
We will heavily make use of a Tauberian Theorem (for $s\to0$, $t\to\infty$, see \cite[Theorem 2 in Chapter XIII.5]{Feller:1971}) which we here quote for the convenience of the reader.
\begin{theorem}\label{thm:Feller_Tauberian}
Let $w:[0,\infty)\to \mathbb{R}$ be a monotone function and assume
its Laplace-Stieltjes transform $\hat w(s) = \int_0^\infty e^{-st} dw(t)$
exists for all $s$ in the right hand complex plane.
Then for $\rho\geq0$
$$
\hat w(s) \sim \frac{C}{s^{\rho}} \quad\mbox{as\ }\ s\to 0
$$
if and only if
$$
w(t) \sim \frac{C}{\Gamma(\rho+1)} t^\rho \quad\mbox{as\ }\ t \to \infty.
$$
\end{theorem}
Note that here $dw$ is actually a measure, which in case of an absolutely
continuous function $w$ can be written as $dw(t)=w'(t)dt$.
\\
Applying Theorem~\ref{thm:Feller_Tauberian} to \eqref{eqn:hath}} 
we get 
\begin{equation}\label{eqn:S_largetime}
h(t)= \frac{B\varphi}{c^2} \sum_{k=1}^N \frac{b_k\lambda^{\beta_k-1}}{\Gamma(1-\alpha_k)} t^{-\alpha_k}  + O(t^{-2\alpha_1}) \mbox{ as }t\to\infty\,.
\end{equation}
Therefore we get an asymptotic formula (like the one in \cite{HatanoNakagawaWangYamamoto:2013} for the subdiffusion equation) for the smallest order 
\begin{equation}\label{eqn:alpha-larget}
\alpha_1 = -\lim_{t\to\infty}\frac{\log h(t)}{\log t}.
\end{equation}
By l'Hospital's rule and actually thereby removing the constant we can instead also use 
\[
\alpha_1 = -\lim_{t\to\infty}\frac{\tfrac{d}{dt}(\log h(t))}{\tfrac{d}{dt}(\log t)} = -\lim_{t\to\infty}\frac{t h'(t)}{h(t)}.
\]
After having determined $\alpha_1$ this way, we can also compute $b_1$ as a limit
\begin{equation}\label{eqn:b1-larget}
b_1 = \lim_{t\to\infty}(h(t) t^{\alpha_1})\, c^2\lambda^{1-\beta_1}\Gamma(1-\alpha_1).
\end{equation}
So we can successively (by the above procedure and subtracting the recovered terms one after another) reconstruct those terms $k$ for which $\alpha_k < 2\alpha_1$, see the remainder term in \eqref{eqn:S_largetime}. However, if there are terms with  $\alpha_k \geq 2\alpha_1$, they seem to get masked by the $O(t^{-2\alpha_1})$ remainder.

The same holds true if we do the excitation by $u_0=0$, $u_1=\varphi$, $f=0$, or $u_0=0$, $u_1=0$, $f=\varphi$, where \eqref{eqn:hath} changes to 
$\hat{h}(s) = \frac{1}{s^2 + \sum_{k=1}^N b_k\lambda^{\beta_k} s^{\alpha_k} + c^2\lambda} B\varphi$.

In order to avoid the restriction $\alpha_k < 2\alpha_1$, we thus refine the expansion of $\hat{h}(s)$ in terms of powers of $s$ and retain the singular ones, that is, those with negative powers, since we are looking at the limiting case $s\to0$.
Using the geometric series formula and the multinomial theorem, with the abbreviations $\Sigma=\sum_{k=1}^N b_k\lambda^{\beta_k} s^{\alpha_k}$, $q=\frac{s^2+\Sigma}{c^2\lambda}$, yields
\[
\begin{aligned}
\hat{h}(s)&=\frac{s+\Sigma/s}{s^2+\Sigma +c^2\lambda} 
= \frac{1}{s}\, \frac{q}{q+1} 
= -\frac{1}{s}\,\sum_{m=1}^\infty (-q)^m
= -\frac{1}{s}\, \sum_{m=1}^\infty (-\tfrac{1}{c^2\lambda})^m\, 
\sum_{\ell=0}^m s^{2(m-\ell)}\, \Sigma^\ell\\
&= -\sum_{m=1}^\infty (-\tfrac{1}{c^2\lambda})^m\, 
\sum_{\ell=0}^m s^{2(m-\ell)-1}\,
\sum_{i_1+\cdots+i_N=\ell}\left({\ell\atop i_1,\cdots,i_N}\right)
\prod_{j=1}^N b_j^{i_j}\, \lambda^{\sum_{j=1}^N\beta_j i_j}
\, s^{\sum_{j=1}^N\alpha_j i_j}.
\end{aligned}
\]
The terms corresponding to $\ell<m$ are obviously $O(s)$, so to extract the singularities it suffices to consider $\ell=m$. The set of indices leading to singular terms is then
\[
I_m=I_m(\alpha_1,\cdots,\alpha_N)=\{\vec{i}=(i_1,\cdots,i_N)\, : \, 
i_1+\cdots+i_N=m, \, \sum_{j=1}^N\alpha_j i_j <1\}\,, \quad
m\leq m_{\max}=[\tfrac{1}{\alpha_1}]
\]
and the singular part of $\hat{h}$ can be written as
\begin{equation}\label{eqn:hsing}
\hat{h}_{\textup{sing}}(s)
= -\sum_{m=1}^{m_{\max}} (-\tfrac{1}{c^2\lambda})^m\, 
\sum_{\vec{i}\in I} \tilde{b}_{m,\vec{i}}
\,\, s^{\sum_{j=1}^N\alpha_j i_j-1}
\mbox{ with } 
\tilde{b}_{\vec{i}} = \left({m\atop i_1,\cdots,i_N}\right)
\prod_{j=1}^N b_j^{i_j}\, \lambda^{\sum_{j=1}^N\beta_j i_j}.
\end{equation}
In case $N=2$ with $i=i_1$, $i_2=m-i$, this reads as
\[
\hat{h}_{\textup{sing}}(s)
= -\sum_{m=1}^{m_{\max}} (-\tfrac{1}{c^2\lambda})^m\, 
\sum_{i=i_{m,\min}}^m \tilde{b}_{m,i}
\,\, s^{\alpha_1 i + \alpha_2(m-i)-1}
\mbox{ with }
\tilde{b}_{m,i} = \left({m\atop i}\right)
b_1^i\, b_2^{m-i}\, \lambda^{\beta_1 i +\beta_2(m-i)}.
\]
since
\[
I_m=\{i\in\{0,\cdots,m\}\, : \, \alpha_1 i +\alpha_2(m-i) <1\}
=\{i_{m,\min},\cdots,m\} \mbox{ with } i_{m,\min}=\left[\frac{m\alpha_2 -1}{\alpha_2-\alpha_1}\right].
\] 
For \eqref{eqn:hsing}, by 
\revision{Theorem~\ref{thm:Feller_Tauberian}} 
\[
\hat{h}_{\textup{sing}}(t)
= -\sum_{m=1}^{m_{\max}} (-\tfrac{1}{c^2\lambda})^m\, 
\sum_{\vec{i}\in I} \frac{\tilde{b}_{m,\vec{i}}}{\Gamma(1-\sum_{j=1}^N\alpha_j i_j)}
\,\, t^{-\sum_{j=1}^N\alpha_j i_j}\,.
\]
which confirms \eqref{eqn:alpha-larget}, \eqref{eqn:b1-larget}.

For the first damping term, we also get an asymptotic formula for the smallest order in Laplace domain -- the small $s$ counterpart to \eqref{eqn:alpha-larget}: Due to
\begin{equation}\label{eqn:asymp-smalls}
\log \hat{h}(s) \approx  \log b_1 +(\alpha_1-1)\log s -  \log(c^2\lambda^{1-\beta_1})\mbox{ as } s\to0
\end{equation}
we have 
\begin{equation}\label{eqn:alpha-smalls}
1-\alpha_1 = -\lim_{s\to0}\frac{\log \hat{h}(s)}{\log s}.
\end{equation}
One might realise this limit by extrapolation: Fitting the values $\frac{\log \hat{h}(s_m)}{\log s_m}$ at $M$ sample points $s_1,\ldots,s_M$ to a regression line or low order polynomial $r(s)$ and setting  $\alpha_1=1+r(0)$.
Next, recover $b_1$ from 
\[
\log b_1 \approx \log \hat{h}(s) + \log(c^2\lambda^{1-\beta_1})-(\alpha_1-1)\log s.
\]
The formula \eqref{eqn:asymp-smalls} also shows why recovering $b_1$ and $\alpha_1$ {\it simultaneously} appears to be so hard: The factor multiplied with $\log b_1$ is unity while the factor multiplied with $1-\alpha_1$ tends to infinity. 

\subsection{Small time asymptotics}
In case $u_0=\varphi$, $u_1=0$, $f=0$ where the Laplace transformed observations are given by 
\begin{equation}\label{eqn:hath_u0}
\hat{h}(s) := \frac{s + \sum_{k=1}^N b_k\lambda^{\beta_k} s^{\alpha_k-1}}{s^2 + \sum_{k=1}^N b_k\lambda^{\beta_k} s^{\alpha_k} + c^2\lambda}
\end{equation}
large $s$ / small $t$ asymptotics cannot be exploited for the following reason.
As $s\to \infty$ then $\hat{h}(s) \to 0$ and in fact $s\hat{h}(s) \to 1$.
More precisely, for
$$d(s) := s\hat{h}(s)-1= \frac{s^2 + \sum_{k=1}^N b_k \lambda^{\beta_k}s^{\alpha_k}}{s^2 + \sum_{k=1}^N b_k \lambda^{\beta_k} s^{\alpha_k} + c^2\lambda} - 1 
= -\frac{c^2\lambda}{s^2 + \sum_{k=1}^N b_k \lambda^{\beta_k} s^{\alpha_k} + c^2\lambda}.
$$
Thus clearly $s^{\gamma} [s\hat{h}(s)-1] \to 0$ for any $\gamma\in[0,2)$.

Not even large $\lambda$ values (or some combined asymptotics as $s\to0$, $\lambda_j\to\infty$) seems to work because then the $c^2\lambda$ term takes over and again masks the terms containing $\alpha_k$.

\medskip
However, in case $u_0=0$, $u_1=\varphi$, $f=0$ the situation is different since then 
\begin{equation}\label{eqn:hath_u1}
\hat{h}(s) := \frac{1}{s^2 + \sum_{k=1}^N b_k\lambda^{\beta_k} s^{\alpha_k} + c^2\lambda}
\end{equation}
and therefore 
\begin{equation}\label{eqn:hhat_large_s}
1-(s^2+c^2\lambda)\hat{h}(s) := 
\underbrace{\frac{s^2}{s^2 + \sum_{k=1}^N b_k\lambda^{\beta_k} s^{\alpha_k} + c^2\lambda}}_{\to 1} 
\Bigl(\sum_{k=1}^N b_k\lambda^{\beta_k} s^{\alpha_k-2}\Bigr)
\end{equation}
which via 
\revision{Theorem~\ref{thm:Feller_Tauberian} (where due to \cite[Theorem 3 in Chapter XIII.5]{Feller:1971} we can replace $s\to0$, $t\to\infty$ by $s\to\infty$, $t\to0$) yields 
}
\[
-c^2\lambda h(t)-h''(t) \sim 
\sum_{k=1}^N \frac{b_k\lambda^{\beta_k}}{\Gamma(2-\alpha_k)} t^{1-\alpha_k}
\mbox{ as }t\to0,
\]
where we have used the fact that $h(0)=0$, $h'(0)=B\varphi$ because of $u_0=0$, $u_1=\varphi$.
\\
This yields the asymptotic formulas 
\[
1-\alpha_N = \lim_{t\to0}\frac{\ln(-c^2\lambda h(t)-h''(t))}{\ln(t)}
= \lim_{t\to0}\left|\frac{t(c^2\lambda h(t)+h''(t))}{c^2\lambda h(t)+h''(t)}\right|
\]
where we have used l'Hospital's rule and 
\[
b_N = -\lim_{t\to0}(c^2\lambda h(t)+h''(t)) t^{\alpha_N-1})\, \lambda^{-\beta_N}\Gamma(2-\alpha_N).
\]
Again, in order to also recover $\alpha_{N-1}$, \ldots, $\alpha_1$, one may also take into account mixed terms in the expansion of \eqref{eqn:hhat_large_s}
analogously to the large $t$ asymptotics case.
In case of two damping terms, this results in  
\begin{equation}\label{eqn:smalltime}
\begin{aligned}
&-c^2\lambda h(t)-h''(t) \\
&=\frac{b_1\lambda^{\beta_1}}{\Gamma(2-\alpha_1)} t^{1-\alpha_1}
+\frac{b_2\lambda^{\beta_2}}{\Gamma(2-\alpha_2)} t^{1-\alpha_2}
+\frac{b_1\lambda^{\beta_1}c^2\lambda}{\Gamma(4-\alpha_1)} t^{3-\alpha_1}
+\frac{b_2\lambda^{\beta_2}c^2\lambda}{\Gamma(4-\alpha_2)} t^{3-\alpha_2}
\\&\quad+\frac{b_1^2\lambda^{2\beta_1}}{\Gamma(4-2\alpha_1)} t^{3-2\alpha_1}
+\frac{b_1b_2\lambda^{\beta_1+\beta_2}}{\Gamma(4-\alpha_1-\alpha_2)} t^{3-\alpha_1-\alpha_2}
+\frac{b_2^2\lambda^{2\beta_2}}{\Gamma(4-2\alpha_2)} t^{3-2\alpha_2}
\\&\quad+ O(t^{5-3\alpha_2})\mbox{\ \  as }\ t\to0 .
\end{aligned}
\end{equation}

%% file: reconstructions.tex
\section{Numerical reconstructions}\label{sec:reconstructions} 

In this section we give some details and show some reconstructions
of the algorithms described in section \ref{sec:recon_meth}.
We assume single mode excitation, set $\beta_i=1$ and
break the reconstructions into three separate paradigms as previously indicated
depending on which time ranges are feasible to be measured.

\input recon_figs.tex
We assume single mode excitation, set $\beta_i=1$ and 
break the reconstructions into three separate paradigms as indicated
in Section~\ref{sec:recon_meth} 
depending on which time ranges are feasible to be measured.

\noindent
\subsection{Full time measurements.}

\begin{figure}[ht]
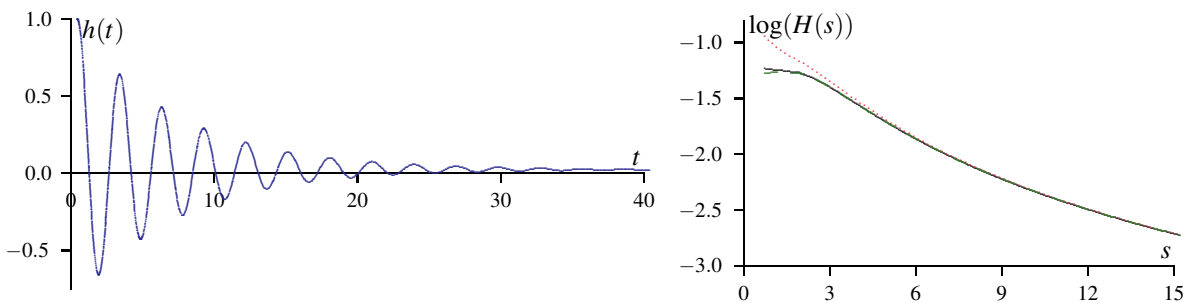

\hbox to \hsize{\hfill\copy\figureone\hfill\copy\figuretwo\hfill}
\small
\caption{\small {\bf Profiles of $h(t)$ and $\hat{h}(s)$, actual and after
iterations 1 and 2.}}
\label{fig:h_and_hath}
\end{figure}

The leftmost graph in Figure \ref{fig:h_and_hath} shows the actual time trace
of $h(t)$ 
over the range $0\leq t\leq 40$
for a damping model with three terms 
$b_i\partial_t^{\alpha_i}$, 
$i=1,2,3$.
Note the damped oscillatory behaviour evident and the range shown
is the one we are labelling as ``full-time'' despite the fact that we
will later look at the important case of 
the very long term
behaviour of the solution where the time values are substantially outside
of this range.  This time range terminology is highly dependent on the
values of critical constants in the model.
These include the size of the domain (which determines the scaling of the
eigenvalues $\{\lambda_n\}$), the wave speed $c$ and of course the strength
of the damping given by the coefficients $\{b_i\}$.
In our reconstructions we will typically take these constants to be of
approximate order unity but note the dependence on these quantities
and the effect that a substantial change would make.

To obtain this $h(t)$ our damping constants are of order 
$\sim 0.1-0.2$ and the combined term 
$\Lambda=c^2\lambda$ 
is taken to be unity.
Changes to the former would modify the time-decrease in $h$ and changes to the 
latter would alter the frequency of the oscillations.

The rightmost graph in Figure~\ref{fig:h_and_hath} shows the logarithm of the
Laplace transform $\hat{h}(s)$ together with the values of $\hat{h}_i(s)$ after
iteration $i$ of the scheme.

This demonstrates two important facts.
First, the very fast convergence of the scheme in the sense of the
convergence of the target Laplace transform function as the parameters
$\{\alpha_i,b_i\}$ are resolved.
The actual $\hat{h}(s)$ is shown by the solid black line, the initial approximation
by the dotted red curve and the first iteration by the dashed green curve
which, at this scale is already almost indistinguishable from the actual $\hat{h}$.
Second, given this it should be quite feasible to obtain reasonably
accurate values the parameters $\{\alpha_i,\,b_i\}$.
In addition, assuming we excite the system with an initial condition
equal to an eigenfunction corresponding to an actual eigenvalue $\lambda$,
we will be able to reconstruct the value of $c^2$ from 
the composite $\Lambda=c^2\lambda$. 

\smallskip
In this situation, to reconstruct the parameters $\{\alpha_i,\,b_i,\, c^2\}$
we work directly with the representation $\hat{h}(s)$ or more exactly with its
logarithm which is a more convenient form for computing the Jacobian
\begin{equation}\label{eqn:log_F(s)}
\log\bigl(\hat{h}(s)\bigr) =
\log\Bigl( s + \sum_1^n b_i s^{\alpha_i-1}\Bigr) -
\log\Bigl( s^2 + \sum_1^n b_i s^{\alpha_i} + c^2\lambda\Bigr).
\end{equation}
Of course, to obtain $\hat{h}(s)$ we must approximate
the integral $\hat{h}(s) = \int_0^\infty e^{-st}h(t)\,dt$ and this indeed does
require full time measurements.
However, we do not actually need the values of the analytic function
$\hat{h}(s)$ for all $s$ for the Newton scheme and so if we only use $s$
values with $s>s_0$ then this allows us to ignore very large time measurements
and indeed we truncated these to $t\leq 40$.

\begin{figure}[h]
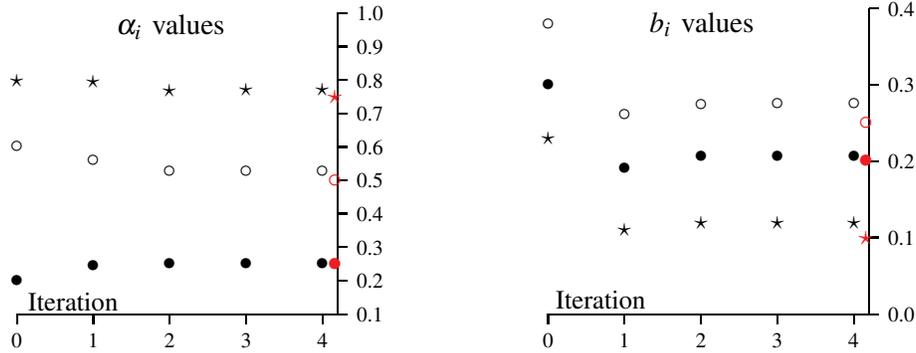

\small
\hbox to \hsize{\hfill\copy\figurefour\hfill\copy\figurethree\hfill}
\caption{\small {\bf reconstructed values of $\{b_i\}$ and $\{\alpha_i\}$.
The symbols in red are the exact values.}}
\label{fig:b_and_alpha}
\end{figure}

The reconstruction of the components in the function $\hat{h}(s)$ is shown
graphically in Figure~\ref{fig:b_and_alpha} and in tabular form in 
Table~\ref{Table:Full_time_trace}.
The exact values were $\alpha=\{0.25, \, 0.5, \, 0.75\}$, $b=\{0.2,\,0.25,\,0.1\}$, $\Lambda=4$. 
Note that we have included the value of the residual here and keeping track of
this value allows the iteration scheme to terminate when saturation has
occurred.
\begin{table}[H]
\centering
\small
\footnotesize
\begin{tabular}{|c|c|c|c|c|c|c|c|c|}
\hline
iter & $\alpha_1$ & $\alpha_2$ & $\alpha_3$ & $b_1$ & $b_2$ & $b_3$ & $\Lambda$ & residual\\
\hline
	0 & 0.3000 & 0.6000 & 0.8000 & 0.3000 & 0.3750 & 0.1500 & 3.500 & {} \\
1 & 0.2448 & 0.5590 & 0.7946 & 0.1907 & 0.2612 & 0.1098 & 4.003 & 0.035032 \\
2 & 0.2506 & 0.5254 & 0.7695 & 0.2057 & 0.2747 & 0.1190 & 4.000 & 0.004371 \\
3 & 0.2492 & 0.5253 & 0.7700 & 0.2060 & 0.2784 & 0.1190 & 4.000 & 0.000075 \\
4 & 0.2491 & 0.5254 & 0.7700 & 0.2060 & 0.2748 & 0.1192 & 4.000 & 0.000000 \\
\hline
\end{tabular}
\small
\footnotesize
\caption{\bf Recovery of  damping terms and unknown $\Lambda$ from full time values.}
 \label{Table:Full_time_trace}
\end{table}

\subsection{Large time measurements.}

Here we are trying to simulate the asymptotic values of the constituent
powers of $t$ occurring in the data function $h(t)$.
This was achieved by using a sample of points between $t_{\rm min}$
and $t_{\rm max}$.

\revision{\subsubsection*{A Newton based approach}}
We apply Newton's method to recover the constants $\{p_i,c_{k,\ell}\}$ in 
\[
\begin{aligned}
h(t) &= 
c_{1,1} t^{-p_1}
+c_{1,2} t^{-p_2}
+c_{1,3} t^{-p_3}
+c_{2,1} t^{-2p_1}
+c_{2,2} t^{-2p_2}
+c_{2,3} t^{-2p_3}
+c_{2,4} t^{-p_1-p_2}
\\&\quad+c_{2,5} t^{-p_2-p_3}
+c_{2,6} t^{-p_3-p_1}
+c_{3,1} t^{-3p_1}
+c_{3,2} t^{-3p_2}
+c_{3,3} t^{-2p_1-p_2}
+c_{3,4} t^{-p_1-2p_2}
\\&\quad+c_{3,5} t^{-2p_1-p_3}
+c_{3,6} t^{-p_1-2p_3}
+c_{3,7} t^{-2p_2-p_3}
+c_{3,8} t^{-p_2-2p_3}
+O(t^{-3p_3})
\end{aligned}
\]
which results from the expansion \eqref{eqn:hsing} in case of three terms,
and then recover $\alpha_i=p_i$ and $b_i=c_{1,i}\Gamma(1-p_i)$  
In the above we have neglected the term $t^{-3\alpha_3}$ as this would not arise
from the Tauberian theorem in the case that the largest power
$\alpha_3 \geq \frac{1}{3}$.
We may also have to exclude other terms such as  
$t^{-2\alpha_1-\alpha_2}$
if
$2\alpha_1 + \alpha_2 \geq 1$.
In practice, during the iteration process, terms should be included or excluded
in the code depending on this criterion: we did so by checking
if the argument passed to the $\Gamma$ function would be negative in which
case the term is deleted from use for that iteration step.

Values for the table shown below were
$t_{\rm min}=5\times 10^4$ and $t_{\rm max}=2\times 10^5$.
As a general rule, terms with small $\alpha$ values can be resolved
with a smaller value of $t_{\rm max}$, but for, say the recovery of
a pair of damping terms with $\alpha_i>0.8$, a larger value 
of $t_{\rm max}$ with commensurate accuracy will be needed.

The case of three damping terms 
$\{b_i\partial^{\alpha_i}_t\}_{i=1}^3$ 
with
$\alpha = \{\frac{1}{4},\,\frac{1}{3},\,\frac{2}{3}\}$ and $b_i=0.1$
is shown in Table~\ref{Table:large_tt} below.
The initial starting values were taken to be between ten and thirty percent
of the actual.
These are shown in the line corresponding to iteration $0$.

\begin{table}[ht]
\centering
\footnotesize
\begin{tabular}{|c|c|c|c|c|c|c|c|}
\hline
iter & $\alpha_1$ & $\alpha_2$ & $\alpha_3$ & $b_1$ & $b_2$ & $b_3$ & residual\\
\hline
0 & 0.2000 & 0.3000 & 0.6000 & 0.130 & 0.080 & 0.110 &\\			
1 & 0.2409 & 0.3572 & 0.6168 & 0.112 & 0.089 & 0.107 & 0.548761\\
2 & 0.2477 & 0.3307 & 0.6531 & 0.102 & 0.091 & 0.104 & 0.104843\\
3 & 0.2499 & 0.3346 & 0.6650 & 0.100 & 0.099 & 0.092 & 0.005124\\
4 & 0.2500 & 0.3332 & 0.6668 & 0.100 & 0.099 & 0.090 & 0.003625\\
5 & 0.2500 & 0.3331 & 0.6658 & 0.100 & 0.100 & 0.090 & 0.001134\\
6 & 0.2500 & 0.3332 & 0.6658 & 0.100 & 0.100 & 0.092 & 0.000578\\
7 & 0.2500 & 0.3332 & 0.6660 & 0.100 & 0.100 & 0.093 & 0.000357\\
8 & 0.2500 & 0.3333 & 0.6665 & 0.100 & 0.100 & 0.094 & 0.000240\\
\hline
\end{tabular}
\small
\caption{Large time values with 3 damping terms}
 \label{Table:large_tt}
\end{table}

There are features here that are typical of such reconstructions.
The reconstruction method resolves the lowest fractional power $\alpha_1$
and its coefficient $b_1$ quickly as this term is the most persistent one
for large times:
essential numerical convergence for $\{\alpha_1,\,b_1\}$ is obtained
by the third iteration.
The next lowest power and coefficient lags behind; here $\{\alpha_2,\,b_2\}$ is already at the stated accuracy by the fifth iteration.
In each case the power is resolved faster and more accurately than its
coefficient.
The third term also illustrates this;  the power is essentially resolved by
iteration 8, but in fact its coefficient $b_3$ is not resolved to the third
decimal place until iteration number 30.
This is seen quite clearly in the singular values of the Jacobian:
the largest singular values correspond to the lowest $\alpha$-values
and the smallest to the coefficients of the largest $\alpha$-powers.

As might be expected, resolving terms whose powers are quite close is
in general more difficult.
This is relatively insignificant for low $\alpha$ values.
For example, with $\alpha_1=0.2$ and $0.22<\alpha_2<0.25$ say, correct
resolution will be obtained although the coefficients will take longer to
resolve than indicated in Table~\ref{Table:large_tt}.
On the other hand if, say $\alpha_1=\frac{1}{4}$ and $\alpha_2=0.85$,
$\alpha_3=0.9$, then with the indicated range of time values used the code
will fail to recover this last pair.
If this is sensed and now only a single second power is requested
this will give a good estimate for $\alpha_2$ but its coefficient will
be overestimated.
Also, as indicated previously, $\alpha$ values close to one require
an extended time measurement range to stay closer to the asymptotic regime of $u(t)$.

\revision{
Note that the starting values were at least ten percent away from the exact values and one cannot expect the convergence radius of Newton's methods to be much larger for a problem exhibiting as high nonlinearity as the one at hand.
Still, let us point to the fact that in the single term case we do not need any starting guesses but obtain very good results directly from the asymptotic formulas \eqref{eqn:alpha-larget}, \eqref{eqn:b1-larget}. This might give the idea of using the asymptotics to construct starting guesses in the multi term case and then apply Newton. However, it is unclear how the asymptotics could yield starting values for terms other than the first one. A (theoretical, as it turns out) possibility for exploiting asympotics in the multi term case is described below. 
}
\Margin{Ref 2 (iv)}

\revision{\subsubsection*{Sucessive sequential use of asymptotic formulas}}
A few words are in order about an approach that from the above discussion
might seem a good or even a better alternative.

Since each damping term 
$b_i \partial^{\alpha_i}_t$ 
contributes a time trace term
with large time behaviour $c_i t^{-\alpha_i}$, it is feasible to take
$T$ sufficiently large so that
$c_i t^{-\alpha_i} <\!\!\!< c_1 t^{-\alpha_1}$ for $t>T$, that is, all but the
smallest damping power is negligible, and this can then be recovered.
In successive steps then we subtract this from the data $h$ to
get $h_1(t) = h(t) - c_1t^{-\alpha_1}$ and now seek to recover the next
lowest $\alpha$ power from the large time values of $h_1(t)$  in a range
$(\delta T,T)$ for $0<\delta<1$.
Then these steps can be repeated until there is no discernible signal remaining
in the sample interval $t_{\rm min},t_{\rm max}$.

This indeed works well under the right circumstances
for recovering two $\alpha$ values but the
coefficients $\{b_i\}$ are less well resolved.
It also requires a delicate splitting of the time interval and gives a
much poorer resolution of the two terms in the case where, say
$\alpha_1=0.2$ and $\alpha_2=0.25$  than that recovered from the Newton scheme.
For the recovery of three damping terms this was in general quite
ineffective.

Every time an $\alpha_i$ has been recovered, the remaining signal is significantly smaller than the previous, leading to an equally significant drop in effective accuracy.
Also, even if we just make a small error in the coefficient $b_i$, the relative error that is caused by this becomes completely dominant for large times. 

In short, this is an elegant and seemingly constructive approach to showing uniqueness for a finite number of damping terms. However, it has limited value from a numerical recovery perspective when used under a wide range of parameter values.

\subsection{Small time measurements.}

In this case we are simulating measurements taken over a very limited initial
time range: in fact we take the measurement interval to be $t\in [0, 0.1)$.
The line of attack is to use the known form of $\hat{h}(s)$ for large values of $s$
and convert the powers of $s$ appearing into powers of $t$ for small times
using the Tauberian theorem.
In case of two damping terms
this gives
\begin{equation*}
\begin{aligned}
-c^2\lambda h_{\tiny\rm{small}}(t)-h_{\tiny\rm{small}}''(t) 
&
=c_{1,1} t^{1-\alpha_1}
+c_{1,2} t^{1-\alpha_2}
+c_{2,1} t^{3-\alpha_1}
+c_{2,2} t^{3-\alpha_2}
\\&\quad
+c_{2,3} t^{3-2\alpha_1}
+c_{2,4} t^{3-\alpha_1-\alpha_2}
+c_{2,5} t^{3-2\alpha_2}
\end{aligned}
\end{equation*}
where each term $c_{k,\ell}$ is computed in terms of
$\{\alpha_i\}$, $\{b_i\}$ and $\lambda$, 
cf. \eqref{eqn:smalltime}. 
							         
The values of $\{\alpha_i,c_{k,\ell}\}$ 
were then computed from the data by a Newton
scheme then finally converted back to the derived values of
$\{\alpha_i\}$ and $\{b_i\}$.

The exact values chosen were
$\alpha = \{0.25,\ 0.2\}$, $b = \{0.1,\ 0.1\}$ and the
initial starting guesses were $\alpha = \{0.3,\ 0.16\}$, $b = \{0.08,\ 0.12\}$.
We show the progression of the iteration in 
Table~\ref{tab:smalltime}. 

\begin{table}[H]
\centering
\small
\footnotesize
\begin{tabular}{|c|c|c|c|c|c|}
\hline
iter & $\alpha_1$ & $\alpha_2$ & $b_1$ & $b_2$ & residual\\
\hline
0 & 0.3000 & 0.1600 & 0.0800 & 0.1200 &  \\
1 & 0.2848 & 0.1673 & 0.0821 & 0.1160 & 0.034308 \\
2 & 0.2632 & 0.1944 & 0.0863 & 0.1142 & 0.024472 \\
3 & 0.2554 & 0.2016 & 0.0862 & 0.1140 & 0.002015 \\
4 & 0.2537 & 0.2033 & 0.0861 & 0.1139 & 0.000444 \\
5 & 0.2534 & 0.2036 & 0.0861 & 0.1139 & 0.000059 \\
6 & 0.2534 & 0.2036 & 0.0861 & 0.1139 & 0.000016 \\
7 & 0.2534 & 0.2036 & 0.0861 & 0.1139 & 0.000010 \\
\hline
\end{tabular}
\small
	\caption{{\small\bf Small time values with 2 damping terms.}
\label{tab:smalltime}}
\end{table}

While theory predicts reconstructibility of an arbitrary number of terms in both cases,
there is a clear difference in ability to reconstruct terms between the small time and the large time asymptotics. 
First of all, the method we described effectively only recovers two terms with small time measurements, as compared to three in the large time. 
The $b_i$ coefficients, which are always harder to obtain than the $\alpha_i$ exponents, are much worse than in the small time than in the large time regime.
This is partly explained by the higher degree of ill-posedness due to the necessity of differentiating the data twice.

%% file: recon_figs.tex
\input colordvi
\input pictex
\font\smallsymbol = cmmi8
\newdimen\xfiglen \newdimen\yfiglen
\xfiglen=2.6 true in
\yfiglen=1.6 true in
\newbox\figurelegendone
\newbox\figurelegendtwo
\newbox\figurelegendthree
\newbox\figureone
\newbox\figuretwo
\newbox\figurethree 
\newbox\figurefour
\newbox\figurefive
\newbox\figuresix
\newbox\figureseven
\newbox\figureeight
\newbox\figurenine
%
\setbox\figurelegendone=\hbox{
\beginpicture
  \setcoordinatesystem units <0.067\xfiglen,0.08\yfiglen> 
  \setplotarea x from 0 to 2, y from 0 to 4
\linethickness=0.6pt
\scriptsize
   \setdashes <3pt>  \putrule from 0 3 to 1.25 3
   \setsolid  
   \Blue{\relax \putrule from 0 2 to 1.25 2 \relax}\relax
   \LimeGreen{\relax \putrule from 0 1 to 1.25 1 \relax}\relax
   \Red{\relax \putrule from 0 0 to 1.25 0 \relax}\relax
  \setsolid
  \put {$u_0$ actual}  [l] at 1.5 3
  \put {$\alpha=0.9$}  [l] at 1.5 2
  \put {$\alpha=0.5$}  [l] at 1.5 1
  \put {$\alpha=0.1$}  [l] at 1.5 0
\endpicture
}
\setbox\figurelegendtwo=\hbox{
\beginpicture
  \setcoordinatesystem units <0.07\xfiglen,0.08\yfiglen> 
  \setplotarea x from 0 to 2, y from 0 to 2
\linethickness=0.6pt
\scriptsize
   \setsolid  
   \Blue{\relax \putrule from 0 2 to 1.25 2 \relax}\relax
   \LimeGreen{\relax \putrule from 0 1 to 1.25 1 \relax}\relax
   \Red{\relax \putrule from 0 0 to 1.25 0 \relax}\relax
  \put {$\alpha=0.9$}  [l] at 1.5 2
  \put {$\alpha=0.5$}  [l] at 1.5 1
  \put {$\alpha=0.1$}  [l] at 1.5 0
\endpicture
}
\setbox\figurelegendthree=\hbox{
\beginpicture
  \setcoordinatesystem units <0.07\xfiglen,0.08\yfiglen> 
  \setplotarea x from 0 to 2, y from 0 to 2
\linethickness=0.6pt
\scriptsize
   \setsolid  
   \Blue{\relax \putrule from 0 2 to 1.25 2 \relax}\relax
   \LimeGreen{\relax \putrule from 0 1 to 1.25 1 \relax}\relax
   \Red{\relax \putrule from 0 0 to 1.25 0 \relax}\relax
  \multiput {\Red{$\bullet$}} at 0.1 0 0.6 0 1.15 0 /
  \multiput {\LimeGreen{$\circ$}}   at 0.1 1 0.6 1 1.15 1 /
  \multiput {\Blue{$\star$}}   at 0.1 2 0.6 2 1.15 2 /
  \put {$\alpha=0.9$}  [l] at 1.5 2
  \put {$\alpha=0.5$}  [l] at 1.5 1
  \put {$\alpha=0.1$}  [l] at 1.5 0
\endpicture
}
\xfiglen=3true in
\yfiglen=0.9true in
\setbox\figureone=\vbox{\hsize=\xfiglen
\beginpicture
  \setcoordinatesystem units <0.025\xfiglen,0.9\yfiglen>  point at 0.0 -0.6
  \setplotarea x from 0 to 40, y from -0.75 to 1
\scriptsize
  \axis bottom shiftedto y=0.0 ticks short numbered from  0 to 40 by 10 /
  \axis left ticks short numbered from -0.5 to 1 by 0.5 /
\footnotesize
\put {$h(t)$} [lt] at 0.6 1
\put {$t$} [rb] at 40 0.05
\setquadratic
\Blue{ 
\input f_full_time.tex
}\relax
\endpicture
}
\xfiglen=2.5true in
\setbox\figuretwo=\vbox{\hsize=\xfiglen
\beginpicture
\small
  \setcoordinatesystem units <0.06\xfiglen,0.65\yfiglen>  point at 0.0 -3
\scriptsize
  \setplotarea x from 0 to 15, y from -3 to -0.8
  \axis bottom shiftedto y=-3 ticks short numbered from 0 to 15 by 3 /
  \axis left ticks short numbered from -3 to -1 by 0.5 /
\footnotesize
\put {$s$} [rb] at 15 -2.93
\put {$\log(H(s))$} [l] at 0.07 -0.85
\setquadratic
\setsolid
\Black{  
\input F_full_s.tex
}\relax
\setdots <2pt>
\Red{\relax   
\input F_full_s_iter1.tex
}\relax
\setdashes
\OliveGreen{\relax   
\input F_full_s_iter2.tex
}\relax
\endpicture
}
%

\xfiglen=0.4true in
\yfiglen=4.0true in
\setbox\figurethree=\vbox{\hsize=\xfiglen   
\beginpicture
  \setcoordinatesystem units <\xfiglen,\yfiglen>  point at 0.0 0.0
  \setplotarea x from 0 to 4.2, y from 0.0 to 0.4
\scriptsize
  \axis bottom ticks short numbered from 0 to 4 by 1 /
  \axis right ticks short numbered from 0.0 to 0.4 by 0.1 /
\footnotesize
\put {Iteration} [lb] at 0.1 0.007
\multiput {$\bullet$} at 0 0.3  1 0.1907  2 0.2057  3 0.2060  4  0.2060 /
\multiput {$\circ$} at 0 0.38  1 0.2612  2 0.2747  3 0.2748  4  0.2748 /
\multiput {$\star$} at 0 0.23  1 0.1098  2 0.1190  3 0.1192  4  0.1192 /
%
\small
\put {\Red{$\bullet$}} at 4.16 0.2
\put {\Red{$\circ$}} at 4.16 0.25
\put {\Red{$\star$}} at 4.16 0.10
\put{$b_i$ values} at 2 0.38
\endpicture
}
\xfiglen=0.4true in
\yfiglen=1.75true in
\setbox\figurefour=\vbox{\hsize=\xfiglen   
\beginpicture
  \setcoordinatesystem units <\xfiglen,\yfiglen>  point at 0.0 0.1
  \setplotarea x from 0 to 4.2, y from 0.1 to 1
\scriptsize
  \axis bottom ticks short numbered from 0 to 4 by 1 /
  \axis right ticks short numbered from 0.1 to 1 by 0.1 /
\footnotesize
\put {Iteration} [lb] at 0.15 0.11
\multiput {$\bullet$} at 0 0.2  1 0.2448  2 0.2506  3 0.2492  4  0.2491 /
\multiput {$\circ$} at 0 0.60  1 0.5590  2 0.5254  3 0.5253  4  0.5253 /
\multiput {$\star$} at 0 0.8  1 0.7946  2 0.7695  3 0.7700  4  0.7700 /
%
\small
\put {\Red{$\bullet$}} at 4.16 0.25
\put {\Red{$\circ$}} at   4.16 0.5
\put {\Red{$\star$}} at   4.16 0.75
\put{$\alpha_i$ values} at 2 0.96
\endpicture
}
\xfiglen=2.5true in
\yfiglen=2true in
\setbox\figurefive=\vbox{\hsize=\xfiglen   
\beginpicture
  \setcoordinatesystem units <10\xfiglen,10\yfiglen>  point at 0.0 0.0
  \setplotarea x from 0 to 0.1, y from 0 to 0.12
\scriptsize
  \axis bottom ticks short numbered from 0 to 0.1 by 0.05 /
  \axis right ticks short numbered from 0 to 0.12 by 0.02 /
\footnotesize
\put {$w(t)$} [lt] at 0.002 0.1
\put {$t$} [rb] at 0.1 0.003
\setquadratic
\Black{ 
\input w_small_time.tex
}\relax
\endpicture
}

%% file: f_full_time.tex
\plot    
         0    1.0000
    0.2000    0.9219
    0.4000    0.7048
    0.6000    0.3927
    0.8000    0.0443
    1.0000   -0.2782
    1.2000   -0.5202
    1.4000   -0.6443
    1.6000   -0.6358
    1.8000   -0.5041
    2.0000   -0.2797
    2.2000   -0.0078
    2.4000    0.2608
    2.6000    0.4787
    2.8000    0.6099
    3.0000    0.6360
    3.2000    0.5580
    3.4000    0.3953
    3.6000    0.1809
    3.8000   -0.0445
    4.0000   -0.2408
    4.2000   -0.3752
    4.4000   -0.4277
    4.6000   -0.3936
    4.8000   -0.2840
    5.0000   -0.1225
    5.2000    0.0595
    5.4000    0.2286
    5.6000    0.3557
    5.8000    0.4206
    6.0000    0.4152
    6.2000    0.3444
    6.4000    0.2241
    6.6000    0.0783
    6.8000   -0.0660
    7.0000   -0.1836
    7.2000   -0.2553
    7.4000   -0.2710
    7.6000   -0.2313
    7.8000   -0.1462
    8.0000   -0.0334
    8.2000    0.0857
    8.4000    0.1897
    8.6000    0.2610
    8.8000    0.2889
    9.0000    0.2708
    9.2000    0.2124
    9.4000    0.1262
    9.6000    0.0289
    9.8000   -0.0618
   10.0000   -0.1303
   10.2000   -0.1659
   10.4000   -0.1641
   10.6000   -0.1273
   10.8000   -0.0642
   11.0000    0.0128
   11.2000    0.0892
   11.4000    0.1516
   11.6000    0.1897
   11.8000    0.1981
   12.0000    0.1770
   12.2000    0.1317
   12.4000    0.0715
   12.6000    0.0078
   12.8000   -0.0480
   13.0000   -0.0866
   13.2000   -0.1021
   13.4000   -0.0931
   13.6000   -0.0626
   13.8000   -0.0172
   14.0000    0.0342
   14.2000    0.0822
   14.4000    0.1186
   14.6000    0.1375
   14.8000    0.1366
   15.0000    0.1170
   15.2000    0.0833
   15.4000    0.0422
   15.6000    0.0013
   15.8000   -0.0323
   16.0000   -0.0530
   16.2000   -0.0579
   16.4000   -0.0471
   16.6000   -0.0234
   16.8000    0.0083
   17.0000    0.0420
   17.2000    0.0715
   17.4000    0.0920
   17.6000    0.1002
   17.8000    0.0955
   18.0000    0.0792
   18.2000    0.0550
   18.4000    0.0275
   18.6000    0.0017
   18.8000   -0.0180
   19.0000   -0.0284
   19.2000   -0.0282
   19.4000   -0.0181
   19.6000   -0.0004
   19.8000    0.0212
   20.0000    0.0429
   20.2000    0.0606
   20.4000    0.0716
   20.6000    0.0742
   20.8000    0.0684
   21.0000    0.0557
   21.2000    0.0388
   21.4000    0.0207
   21.6000    0.0048
   21.8000   -0.0063
   22.0000   -0.0109
   22.2000   -0.0086
   22.4000   -0.0002
   22.6000    0.0125
   22.8000    0.0270
   23.0000    0.0406
   23.2000    0.0510
   23.4000    0.0564
   23.6000    0.0563
   23.8000    0.0509
   24.0000    0.0414
   24.2000    0.0298
   24.4000    0.0181
   24.6000    0.0085
   24.8000    0.0025
   25.0000    0.0010
   25.2000    0.0039
   25.4000    0.0104
   25.6000    0.0193
   25.8000    0.0288
   26.0000    0.0372
   26.2000    0.0431
   26.4000    0.0455
   26.6000    0.0442
   26.8000    0.0396
   27.0000    0.0328
   27.2000    0.0250
   27.4000    0.0176
   27.6000    0.0120
   27.8000    0.0090
   28.0000    0.0089
   28.2000    0.0116
   28.4000    0.0165
   28.6000    0.0226
   28.8000    0.0287
   29.0000    0.0338
   29.2000    0.0369
   29.4000    0.0377
   29.6000    0.0361
   29.8000    0.0325
   30.0000    0.0277
   30.2000    0.0226
   30.4000    0.0180
   30.6000    0.0148
   30.8000    0.0134
   31.0000    0.0140
   31.2000    0.0163
   31.4000    0.0198
   31.6000    0.0239
   31.8000    0.0277
   32.0000    0.0307
   32.2000    0.0322
   32.4000    0.0322
   32.6000    0.0307
   32.8000    0.0280
   33.0000    0.0247
   33.2000    0.0214
   33.4000    0.0186
   33.6000    0.0168
   33.8000    0.0163
   34.0000    0.0171
   34.2000    0.0189
   34.4000    0.0214
   34.6000    0.0240
   34.8000    0.0264
   35.0000    0.0281
   35.2000    0.0287
   35.4000    0.0284
   35.6000    0.0271
   35.8000    0.0252
   36.0000    0.0229
   36.2000    0.0208
   36.4000    0.0192
   36.6000    0.0183
   36.8000    0.0182
   37.0000    0.0189
   37.2000    0.0203
   37.4000    0.0220
   37.6000    0.0237
   37.8000    0.0251
   38.0000    0.0260
   38.2000    0.0262
   38.4000    0.0257
   38.6000    0.0247
   38.8000    0.0233
   39.0000    0.0219
   39.2000    0.0206
   39.4000    0.0196
   39.6000    0.0192
   39.8000    0.0193
   40.0000    0.0199
/\relax

%% file: F_full_s.tex
\plot    
    0.5000   -1.2306
    1.1042   -1.2509
    1.7083   -1.2687
    2.3125   -1.3327
    2.9167   -1.4209
    3.5208   -1.5178
    4.1250   -1.6153
    4.7292   -1.7098
    5.3333   -1.7997
    5.9375   -1.8846
    6.5417   -1.9645
    7.1458   -2.0397
    7.7500   -2.1104
    8.3542   -2.1771
    8.9583   -2.2401
    9.5625   -2.2998
   10.1667   -2.3564
   10.7708   -2.4102
   11.3750   -2.4614
   11.9792   -2.5103
   12.5833   -2.5570
   13.1875   -2.6017
   13.7917   -2.6446
   14.3958   -2.6858
   15.0000   -2.7255
/\relax

%% file: F_full_s_iter1.tex
\plot    
    0.5000   -0.9432
    1.1042   -1.0771
    1.7083   -1.1620
    2.3125   -1.2635
    2.9167   -1.3737
    3.5208   -1.4840
    4.1250   -1.5901
    4.7292   -1.6905
    5.3333   -1.7844
    5.9375   -1.8723
    6.5417   -1.9544
    7.1458   -2.0312
    7.7500   -2.1032
    8.3542   -2.1709
    8.9583   -2.2348
    9.5625   -2.2951
   10.1667   -2.3522
   10.7708   -2.4065
   11.3750   -2.4581
   11.9792   -2.5073
   12.5833   -2.5543
   13.1875   -2.5993
   13.7917   -2.6424
   14.3958   -2.6838
   15.0000   -2.7236
/\relax

%% file: F_full_s_iter2.tex
\plot    
    0.5000   -1.2731
    1.1042   -1.2661
    1.7083   -1.2750
    2.3125   -1.3356
    2.9167   -1.4225
    3.5208   -1.5187
    4.1250   -1.6159
    4.7292   -1.7102
    5.3333   -1.8000
    5.9375   -1.8848
    6.5417   -1.9646
    7.1458   -2.0398
    7.7500   -2.1105
    8.3542   -2.1772
    8.9583   -2.2402
    9.5625   -2.2998
   10.1667   -2.3564
   10.7708   -2.4102
   11.3750   -2.4614
   11.9792   -2.5103
   12.5833   -2.5570
   13.1875   -2.6017
   13.7917   -2.6446
   14.3958   -2.6858
   15.0000   -2.7255
/\relax

%% file: w_small_time.tex
\plot    
       0         0
    0.0025    0.0025
    0.0050    0.0050
    0.0075    0.0075
    0.0100    0.0100
    0.0125    0.0125
    0.0150    0.0150
    0.0175    0.0175
    0.0200    0.0200
    0.0225    0.0225
    0.0250    0.0250
    0.0275    0.0275
    0.0300    0.0300
    0.0325    0.0325
    0.0350    0.0350
    0.0375    0.0375
    0.0400    0.0399
    0.0425    0.0424
    0.0450    0.0449
    0.0475    0.0474
    0.0500    0.0499
    0.0525    0.0524
    0.0550    0.0549
    0.0575    0.0574
    0.0600    0.0599
    0.0625    0.0623
    0.0650    0.0648
    0.0675    0.0673
    0.0700    0.0698
    0.0725    0.0723
    0.0750    0.0748
    0.0775    0.0772
    0.0800    0.0797
    0.0825    0.0822
    0.0850    0.0847
    0.0875    0.0871
    0.0900    0.0896
    0.0925    0.0921
    0.0950    0.0945
    0.0975    0.0970
    0.1000    0.0995
/\relax